\documentclass[11pt]{article}
\usepackage[left=2.5cm, right=2.5cm, top=2.5cm, bottom=2.5cm]{geometry}
\usepackage{amsmath,amsfonts,amssymb,amsthm,mathtools}
\usepackage{graphicx}
\usepackage{cite}

\newtheorem{prop}{Proposition}

\newtheorem{property}{Property}

\usepackage{color}
\usepackage{bm}

\usepackage{fancybox}

\begin{document}

\title{Drawn Games and Tied Pools }


\author{R. Edwards\\Department of Mathematics and Statistics, University of Victoria, \\
PO Box 1700, STN CSC, Victoria, BC, Canada, V8W 2Y2\\edwards@uvic.ca}

\maketitle
\begin{abstract}

We consider the problem of estimating `preference' or `strength' parameters in three-way comparison experiments, each composed of a series of paired comparisons, but where only the single `preferred' or `strongest' candidate is known in each trial. Such experiments arise in psychology and market research, but here we use chess competitions as the prototypical context, in particular a series of `pools' between three players that occurred in 1821. The possibilities of tied pools, redundant and therefore unplayed games, and drawn games must all be considered. This leads us to reconsider previous models for estimating strength parameters when drawn games are a possible result. In particular, Davidson's method for ties has been questioned, and we propose an alternative. We argue that the most correct use of this method is to estimate strength parameters first, and then fix these to estimate a draw-propensity parameter, rather than estimating all parameters simultaneously, as Davidson does. This results in a model that is consistent with, and provides more context for, a simple method for handling draws proposed by Glickman. Finally, in pools with incomplete information, the number of drawn games can be estimated by adopting a draw-propensity parameter from related data with more complete information.

 \end{abstract}

\section*{Keywords:} Paired comparisons with ties, Three-way comparisons, Maximum-likelihood estimation, Bradley-Terry model, Chess

\section{Introduction}

Consider the problem of evaluating a `preference' value for each of three images or items by having a number of judges select their favourite by means of successive paired comparisons. While all three items could be presented to a judge simultaneously, one might expect that in the judge's mind a series of pairwise comparisons occurs in any case, and it could be preferrable to make sure that the judges are clear about their own preferences by forcing them to make the comparisons pairwise. 

More specifically, the procedure envisaged is as follows. For each judge, two items are randomly selected for the initial comparison, or equivalently, one of the three items is selected to be left out in the initial comparison. If the three items are labelled $B$, $C$, and $D$, then suppose item $D$ is left out at first. If the judge prefers $B$ to $C$, then the next comparison is $B$ to $D$. At that point if the judge also prefers $B$ to $D$, then the third possible comparison need not be made, since $B$ is preferred to both of the others. However, if $D$ is preferred to $B$, then the final comparison, between $C$ and $D$, is made. If $D$ is preferred to $C$, then $D$ is selected as the favourite. If now $C$ is preferred to $D$, then the preferences are circular ($B$ over $C$, $C$ over $D$, and $D$ over $B$), and the entire procedure must be carried out again. Additional adjustments must be made if a judge at some point finds a pair of items equally preferable. One approach would be to simply present the same two items repeatedly until a decision can be made (as might be done, for example, if a comparison was requested after showing two images to the judge very briefly, for perhaps a few seconds, to elicit only a first impression, not a more considered one). If all the details of the individual pairwise comparisons are known then standard theory for paired comparison experiments can be applied, like the Bradley-Terry model~\cite{Bradley1952}, or to account for ties, the Davidson model~\cite{Davidson1970}. However, if only the favourite of each judge is recorded, then that theory is not sufficient in itself to deal with the three-way comparisons.

The same problem arises in the context of a chess competition involving three players, in the form of pools.  In each round, the three players each put their bet into a pool, and then play each other in pairs successively until one player beats both of the others, and then wins the pool. Exactly this type of event was played by three of the best players of England and France in 1821: Louis Charles Mah\'e de la Bourdonnais, John Cochrane, and Alexandre Louis Honor\'e Lebreton Deschapelles. Using the first letter of each surname, we can label them $B$, $C$, and $D$, respectively\footnote{I use the most common form of the French players' names, suggesting initials B, C, and D, though sometimes variant spellings appear, such as De la Bourdonnais, de Labourdonnais, and Des Chapelles. One complication is that Deschapelles gave the other two players a handicap (of a pawn and two moves to start), which means that he was playing at reduced strength, but we can take his `playing strength' parameter to reflect the weakened Deschapelles, not the full-strength (un-handicapped) Deschapelles.}. The format of the event as a series of pools (`{\em poules}' in French) is described explicitly by Saint-Amant~\cite[p.313]{SaintAmant1844}. The most complete information we have on the results of this event, by Staunton~\cite[p.25]{Staunton1841v1}, is given in Table~\ref{StauntonData}.
\begin{table}
\centering{
\begin{tabular}{|l|r|} \hline
Number of pools won by $B$ & 14 \\ 
Number of pools won by $C$ & 1 \\
Number of pools won by $D$ & 6 \\ \hline \hline
Total number of pools won & 21 \\ \hline
\end{tabular}
} 
\caption{Results of chess pools played at Paris 1821\label{StauntonData}}
\end{table}
A later account~\cite[p.250]{Walker1843} says that that each player bet 1 napoleon in each round. Nowhere are details given of how the play proceeded or who won which individual games. However, it is likely that the player to be initially left out was either selected randomly in each round, or rotated through the rounds (the fact that the number of rounds was a multiple of 3 suggests the latter). It is also likely that redundant games were not played, since if the winner of the pool was determined by the first two pairings, there would be no need for the third to be played. And if the players each won one game in a pool (a tied pool), it would have to be replayed, as no winner was yet determined. Finally, drawn games were similarly considered inconclusive in that era and were simply replayed (for example, in the matches between de la Bourdonnais and McDonnell in 1834~\cite[p.21]{Utterberg2005}, and the Staunton-Harrwitz matches of 1846~\cite[p.60]{Sergeant1934}). Henceforth, we shall distinguish between a pool and a round, where a round includes as many pools as were necessary for someone to win the bet. Thus, there can be tied pools but not tied rounds.

The problem in either of the situations above is to determine overall `preference values' or `player strengths' based on the results of the three-way comparisons, when only the favourite of the three items is known for each judge, or only the winner of each pool is known, and not the details of the individual pairwise comparisons. In the rest of this paper, we will adopt the language of the chess competition, but it should be remembered that the theory applies more generally.

David~\cite{David1988} summarizes several methods for dealing with three-way comparisons, the most intuitive of which is the one by Pendergrass and Bradley~\cite{Pendergrass1960}, which treats three-way comparisons as a sequence of paired comparisons. We follow their general approach here, but develop our ideas from scratch, starting in Section~\ref{sec:pools}, and returning to the 3-way problem in Section~\ref{sec:draws}, after a lengthy assessment of the question of handling drawn games.

David also provides a good summary of methods for estimating preferences or strengths in paired comparison experiments, notably the Bradley-Terry method~\cite{Bradley1952} among linear models, and methods for including ties or draws as possible outcomes of a comparison, notably the method of Davidson~\cite{Davidson1970}, as well as that of Rao and Kupper~\cite{Rao1967}, both of which include a draw-propensity parameter. Amongst more recent theoretical developments in which draws are considered, are those of Glickman and Hankin. Glickman's paper of 1999~\cite{Glickman1999} simply treats a draw as half a win (scoring $0.5$ instead of $0$ or $1$, as is normally done in actual chess competitions), thus avoiding the need for a draw-propensity parameter. Glickman's 2025 paper~\cite{Glickman2025} and Hankin~\cite{Hankin2020} deal also with order effects, which are known to create a small imbalance in chess, but the data presented by Glickman's Figure 2 does not make a convincing case for a significant order effect, and we do not consider it here. In the chess context, in any case, there is usually an attempt to balance the number of times each player has the first move, so that any imbalance due to order is cancelled out. Hankin handles draws by a novel approach of introducing a hypothetical third entity, a `draw monster,' who wins when a game between two players is drawn, designed to focus on the possible collusion that may have occurred among Soviet chess players to draw games against each other and fight for wins against the non-Soviets. Glickman~\cite{Glickman2025} is interested in strength-dependent draw propensity, and modifies Davidson's approach to allow for this. Neither of these issues are our main concern, but we do need to include draw propensity in order to estimate numbers of draws between players when that information is not known, so we start from Davidson's method, but propose an alternative to improve on two aspects of the method that seem flawed, or at least undesirable.

In Section~\ref{sec:pce} we present briefly the Bradley-Terry model and its generalization by Davidson to allow for draws. In Section~\ref{sec:alternative}, we propose an Alternative model that makes a natural-looking correction to an assumption of Davidson's model that was highlighted by David as questionable. We outline the maximum likelihood estimation that results from Davidson's and the Alternative model in Section~\ref{sec:max_likelihood}, but show in Section~\ref{sec:consistency} that the Alternative model does not necessarily rank player strengths in the same order as their scores in balanced designs (tournaments), a desirable property that we call {\em Consistency}, satisfied by Davidson's model. In Section~\ref{sec:reassessment}, we step back and reassess the estimation procedure followed by Davidson, and so far copied in the Alternative model, and conclude that it is more sensible to separate the estimation of strength parameters from estimation of draw propensities; the latter should not influence the former. This leads to four candidate models, the original two, and versions of each of them constrained by the separation of estimation steps. We apply these four models to some illustrative small test cases in Section~\ref{sec:tests}, leading us to conclude that the Constrained Alternative model is most appropriate, satisfying the Consistency property, and correcting the flawed assumptions of the other three models. After integrating draws into our 3-way analysis in Section~\ref{sec:draws}, we finish by carrying out our full procedure on the historical example in Section~\ref{sec:Paris}.

\section{Incomplete information on pools} \label{sec:pools}

\subsection{Without draws}
We first assume that there are no drawn games. We will introduce this possibility later. For now, we wish to deal with tied pools and redundant games.

Let $p_{BC}$ be the probability that $B$ beats $C$ in a game, and similarly for other pairs. Since we do not allow draws for now, $p_{BC}+p_{CB}=1$, {\em etc.} The probability that player $B$ wins a pool is 
$$q_B = p_{BC}p_{BD}\,.$$ 
This holds whether or not a game between $C$ and $D$ was played. Similarly, 
$$q_C=p_{CB}p_{CD}=(1-p_{BC})p_{CD}\,,\mbox{ and}$$ 
$$q_D=p_{DB}p_{DC}=(1-p_{BD})(1-p_{CD})\,.$$
The probability that no player wins a pool, {\em i.e.}, the pool is tied, is
$$p_0=p_{BC}p_{CD}p_{DB}+p_{BD}p_{DC}p_{CB}=1-(q_B+q_C+q_D)\,.$$
Finally, given the possibility of tied pools that have to be repeated, possibly many times until a winner emerges, the probability that $B$ eventually wins a round is
$$p_B=p_{BC}p_{BD}+p_0p_{BC}p_{BD}+p_0^2p_{BC}p_{BD}+\ldots = p_{BC}p_{BD}\sum_{k=0}^{\infty}p_0^k=\frac{p_{BC}p_{BD}}{1-p_0}=\frac{q_B}{1-p_0}.$$

Now, the Bradley-Terry model for paired comparisons (presented more formally in Section~\ref{sec:pce}) assumes that there exist strength parameters for each player, $\pi_B$, $\pi_C$, and $\pi_D$, such that 
$$p_{BC}=\frac{\pi_B}{\pi_B+\pi_C}, \quad p_{BD}=\frac{\pi_B}{\pi_B+\pi_D}, \quad p_{CD}=\frac{\pi_C}{\pi_C+\pi_D},$$
and, since $(\pi_B,\pi_C,\pi_D)$ is clearly defined only up to a linear scaling, we can normalize by setting 
$$\pi_B+\pi_C+\pi_D=1.$$

Our problem is, given proportions of won rounds, $(p_B,p_C,p_D)=\frac{(q_B,q_C,q_D)}{1-p_0}=\frac{(q_B,q_C,q_D)}{q_B+q_C+q_D}$, to find estimates (expected values) of $(\pi_B,\pi_C,\pi_D)$.

In order to do this, consider ratios of the probabilities:
$$\frac{p_B}{p_C}=\frac{q_B}{q_C}=\frac{p_{BC}p_{BD}}{p_{CB}p_{CD}}=\frac{\pi_B}{(\pi_B+\pi_C)}\frac{\pi_B}{(\pi_B+\pi_D)}\frac{(\pi_B+\pi_C)}{\pi_C}\frac{(\pi_C+\pi_D)}{\pi_C}=\frac{\pi_B^2(\pi_C+\pi_D)}{\pi_C^2(\pi_B+\pi_D)}=\frac{\pi_B^2(1-\pi_B)}{\pi_C^2(1-\pi_C)}\,,$$
and letting $f(x)=x^2(1-x)$, we have 
$$\frac{p_B}{p_C}=\frac{f(\pi_B)}{f(\pi_C)},\quad \mbox{ and similarly, } \quad\frac{p_B}{p_D}=\frac{f(\pi_B)}{f(\pi_D)}, \quad \frac{p_C}{p_D}=\frac{f(\pi_C)}{f(\pi_D)}\,,$$
where any one of these can clearly be derived from the other two, so there is a dependency between these three equations. Note that this is similar to the model for triple comparisons suggested by Pendergrass and Bradley~\cite{Pendergrass1960}, which treats triple comparisons as composed of a set of paired comparisons. Given values for $(p_B,p_C,p_D)$, we can obtain $(\pi_B,\pi_C,\pi_D)$ by solving simultaneously the three equations:
\begin{eqnarray*}
p_Cf(\pi_B) - p_Bf(\pi_C) &=& 0\\
p_Cf(\pi_D) - p_Df(\pi_C) &=& 0 \\
\pi_B+\pi_C+\pi_D &=& 1
\end{eqnarray*}
These are nonlinear, but one way is to use Newton's method to find a root of 
$$ F(\pi_B,\pi_C,\pi_D)=[p_Cf(\pi_B) - p_Bf(\pi_C)]^2 + [p_Cf(\pi_D) - p_Df(\pi_C) ]^2 + [1-(\pi_B+\pi_C+\pi_D)]^2\,.$$
Letting $\bm{\pi}=(\pi_B,\pi_C,\pi_D)^\top$, Newton's method iterates
$$ \bm{\pi}_{n+1}=\bm{\pi}_n-{\bf J}^{-1}{\bf F}(\bm{\pi}) $$
where
$$ {\bf F}(\bm{\pi}) =\left[\begin{array}{c} p_Cf(\pi_B) - p_Bf(\pi_C) \\ p_Cf(\pi_D) - p_Df(\pi_C) \\ 1-(\pi_B+\pi_C+\pi_D) \end{array}\right]\quad\mbox{and}\quad {\bf J}= \left[\begin{array}{ccc} 
p_Cf^\prime(\pi_B) & -p_Bf^\prime(\pi_C) & 0 \\
0 & -p_Df^\prime(\pi_C) & p_Cf^\prime(\pi_D) \\
-1 & -1 & -1 \end{array}\right]\,, $$
with $f^\prime(x)=x(2-3x)$. One may start, for example with an initial guess of $\bm{\pi}_0=(p_B,p_C,p_D)^\top$. 

In the historical example, we have proportions of won rounds $(p_B,p_C,p_D)=\left(\frac{14}{21},\frac{1}{21},\frac{6}{21}\right)$, and iteration of Newton's method converges to the solution 
$$\bm{\pi}\approx (0.5971497, 0.1072050, 0.2956435 )^\top\,.$$ 
We can use these values to compute estimates of all of the other quantities. Our estimate of the probabilities of wins in individual games by one player over another are: $p_{BC}\approx 0.8477968,\, p_{BD}\approx 0.6688542,\, p_{CD}\approx 0.2661163$ and our estimates of the probabilites that one player wins a pool are $q_B\approx 0.5670525,\, q_C\approx 0.04050375, \, q_D\approx 0.2430225$. 
The expected proportion of tied pools is $p_0\approx 0.1494213$, so that the expected number of pools played (including replayed tied pools) is $\frac{21}{1-p_0}\approx 24.68907$. 

Finally, we can estimate the number of redundant (and therefore unplayed) games as follows. A game between $C$ and $D$ is not played if $D$ gets a bye (sits out) first, and $B$ beats both (probability $\frac 13 p_{BC}p_{BD}=\frac 13 q_B$) or if $C$ gets a bye first, and $B$ beats both (same probability), so the game between $C$ and $D$ is not played with probability $\frac 23 q_B$. Thus, the expected number of games between $C$ and $D$ overall is 
$$n_{CD}=\left(1-\frac 23 q_B\right)\frac{N}{1-p_0}$$
where $N$ is the total number of rounds ({\em i.e.}, not counting tied pools). Similarly, $n_{BC}=\left(1-\frac 23 q_D\right)\frac{N}{1-p_0}$, and $n_{BD}=\left(1-\frac 23 q_C\right)\frac{N}{1-p_0}$. And the expected score for player $C$ in games against $D$ is
$$s_{CD}=n_{CD}p_{CD}$$
and similarly for the other pairs.

For our historical example, we get $n_{BC}\approx 20.68907, \, n_{BD}\approx 24.02241, \, n_{CD}\approx 15.35574 $ and $s_{BC}\approx 17.54013, \, s_{BD}\approx 16.06749, \, s_{CD}\approx 4.08641$.

\subsection{With draws}
But what if we recognize the possibility of draws? Then $p_{BC}+p_{CB}+d_{BC}=1$, where $d_{BC}=d_{CB}$ is the probability of a draw between $B$ an $C$. Now, if drawn games are simply replayed until one of the two players wins, then the probability that $B$ wins the `encounter' (no longer necessarily a single game) with $C$ is 
\begin{equation} e_{BC}=p_{BC}+d_{BC}p_{BC}+d_{BC}^2p_{BC}+\ldots = p_{BC}\sum_{k=0}^{\infty}d_{BC}^k=\frac{p_{BC}}{1-d_{BC}}\,,
\label{encounterbygame}
\end{equation}
and similarly for the other encounters. 
Then $e_{BC}+e_{CB}=\frac{p_{BC}+p_{CB}}{1-d_{BC}}=1$, and we can redefine $q_B=e_{BC}e_{BD}$, {\em etc.}, and $p_0=1-(q_B+q_C+q_D)$ as before, and $p_B=\frac{q_B}{1-p_0}$, {\em etc.}, the probabilities of each player to win a round. But then, we need to decide how to define strength parameters. Also, the estimates labelled $n_{BC}$ and $s_{BC}$, {\em etc.}, above are no longer counts and scores of games, but of encounters, and to avoid confusion, we should label them differently, as $s^e_{BC}$ and $n^e_{BC}$, {\em etc.}, for example. For later use, we list these values here:
\begin{equation} n^e_{BC}\approx 20.68907, \, n^e_{BD}\approx 24.02241, \, n^e_{CD}\approx 15.35574 \label{eq:ne}
\end{equation}
and 
\begin{equation} s^e_{BC}\approx 17.54013, \, s^e_{BD}\approx 16.06749, \, s^e_{CD}\approx 4.08641\,. \label{eq:se}
\end{equation}

One possibility might be to ignore the draws and base the strength parameters on probabilities of an eventual win of an encounter:
$$e_{BC}=\frac{\pi_B}{\pi_B+\pi_C}\,,\quad\mbox{and}\quad e_{CB}=\frac{\pi_C}{\pi_B+\pi_C}\,.$$
Then $\frac{\pi_B}{\pi_C}=\frac{e_{BC}}{e_{CB}}=\frac{p_{BC}}{p_{CB}}$. 
This has the advantage that the previous theory without draws still applies, and so our estimates of player strength parameters will be the same as in that case. 

However, this assumption ignores the natural idea that many draws between a pair of players indicates similar playing strength. In the context of our sequence of pools, if the score of won games between a pair of players is $15$ to $4$, that suggests a significant difference in playing strength, but if there were a large number of drawn games that preceded each win, that would lessen our sense of the difference in strength, since most of the time neither player could beat the other. So we turn now to a consideration of the modelling of draws in paired comparison experiments in general, and we return to the issue of estimating player strengths from data on the results of pools in Section~\ref{sec:draws}.

\section{Paired comparison experiments, without and with draws}
\label{sec:pce}

We are concerned with paired comparison experiments \cite{Bradley1952,David1988} and methods for handling ties or draws in paired comparison experiments \cite{Davidson1970,Rao1967}. Let $p_{ij}$ be the probability that player $i$ beats player $j$, where $i$ and $j$ range from $1$ to $t$ (David~\cite{David1988} uses $t$ to reprsent the number of `treatments', here the number of players). 

\subsection{The Bradley-Terry model}
If draws are not permitted (or not accounted for) then \begin{equation}
p_{ij}+p_{ji}=1\,.\label{pp}
\end{equation}
The Bradley-Terry model (without draws) supposes that there are player `strength' parameters, $\pi_i$ for player $i$ ($i=1,\ldots,t$), satisfying 
\begin{property}\label{Luce_Condition}
 \begin{equation}\frac{\pi_i}{\pi_j}=\frac{p_{ij}}{p_{ji}}\,,\label{Luce}\end{equation}   
\end{property}
Property~\ref{Luce_Condition} guarantees that
\begin{equation}\frac{p_{ij}}{p_{ji}}\frac{p_{jk}}{p_{kj}}=\frac{p_{ik}}{p_{ki}}\,.\label{bt_davidson}
\end{equation}
Equations~\eqref{pp} and \eqref{Luce} together imply that 
\begin{equation}
    p_{ij}=\frac{\pi_i}{\pi_i+\pi_j}\,.
\end{equation}
Note that this only defines $\bm{\pi}=(\pi_1,\ldots,\pi_t)^\top$ up to an arbitrary scaling. In practice this is often fixed using $\sum_i \pi_i=1$.
Chess ratings usually use a transformation of these strength parameters such that $r_i=r_0+400\log_{10}(\pi_i)$ (where $r_0$ is an offset that depends on the scaling of $\bm{\pi}$).


\subsection{The Davidson model}
Now, let $d_{ij}$ be the probability of a draw between player $i$ and player $j$ (so $d_{ji}=d_{ij}$). Thus, we now have
\begin{equation}p_{ij}+p_{ji}+d_{ij}=1\,,\label{ppd}
\end{equation}
Davidson's model~\cite{Davidson1970} supposes that 
\begin{equation}d_{ij}=\nu \sqrt{p_{ij}p_{ji}}\label{DavidsonDraw}
\end{equation}
with parameter $\nu\ge 0$, and that Property~\ref{Luce_Condition} still holds. 
Equations~\eqref{ppd} and \eqref{DavidsonDraw} together 
imply that
$$p_{ij}+p_{ji}+\nu\sqrt{p_{ij}p_{ji}}=1$$
and dividing through by $p_{ij}$ and rearranging using Equation~\eqref{Luce} gives
\begin{equation}p_{ij}=\frac{\pi_i}{\pi_i+\pi_j+\nu\sqrt{\pi_i\pi_j}}\,, \quad\mbox{and}\quad d_{ij}=\frac{\nu\sqrt{\pi_i\pi_j}}{\pi_i+\pi_j+\nu\sqrt{\pi_i\pi_j}}\,.\end{equation}
When $\nu=0$, the Davidson model reduces to the Bradley-Terry model.

The parameters $\bm{\pi}$ and $\nu$ are estimated from data in the form of game results. Let $w_{ij}$ be the number of games won by player $i$ against player $j$, and let $t_{ij}$ be the number of draws that occurred between them. Clearly $t_{ij}=t_{ji}$. Let $s_{ij}=w_{ij}+\frac 12t_{ij}$, the score of player $i$ in games against player $j$ counting every draw as half a win, and $s_i=\sum_j s_{ij}$, the total score of player $i$ against all opponents. We can denote the total number of games between a pair as $n_{ij}=s_{ij}+s_{ji}=w_{ij}+w_{ji}+t_{ij}$, and the total number of games played by player $i$ as $n_i=\sum_j n_{ij}$. We assume throughout, naturally, that $n_{ii}=0$ for all $i$ (and thus, $w_{ii}=0$, $t_{ii}=0$, and $s_{ii}=0$). Finally, we will have occasion to use $T=\sum_{i<j}t_{ij}=\frac 12 \sum t_{ij}$, the total number of draws between all pairs of players, and
$$N=\frac 12 \sum_i n_i=\sum_{i<j}n_{ij}=\sum_{i<j}(w_{ij}+w_{ji}+t_{ij})=\sum_{i,j}w_{ij}+\sum_{i<j}t_{ij}=\sum_{i,j}s_{ij}=\sum_i s_i\,,$$
the total number of games.

One of the desirable (perhaps essential) properties of a model is that player strengths are ranked according to scores in balanced designs (all $n_{ij}$ equal):
\begin{property}\label{consistency_property}
Consistency: In a balanced design (when all $n_{ij}\equiv n$, for $i\ne j$), $\pi_i-\pi_j$ and $s_i-s_j$ have the same sign, or equivalently, $\log\left(\frac{\pi_i}{\pi_j}\right)$ and $\log\left(\frac{s_i}{s_j}\right)$ have the same sign.
\end{property}
We will see below in Section~\ref{sec:max_likelihood} that Davidson's model satisfies Consistency.

David~\cite[p.140]{David1988}, however, questions Davidson's assumption of Property~(\ref{Luce_Condition}), which implies that the probabilities of draws have no effect on the rating gap between a single pair of players. Considering two players in isolation, we will obtain the estimate $\frac{\pi_1}{\pi_2}=\frac{w_{12}}{w_{21}}$ no matter what $t_{12}$ is (though the estimate of $\nu$ will of course depend on $t_{12}$). Conversely, the same value of $\frac{s_1}{s_2}$ involving different numbers of draws can yield different estimates of the strength parameter ratio.
For example,  if $t_{12}=0$ and $s_1=w_{12}=48$, while $s_2=w_{21}=24$, then the strength ratio $\frac{\pi_1}{\pi_2}$ is estimated to be $2$ (and $\nu=0$), while if $w_{12}=32, w_{21}=8$ and $t_{12}=32$ so that the scores are still $s_1=48$ and $s_2=24$ with $T=32$ draws, then we estimate $\nu=\frac{T}{\sqrt{w_{12}w_{21}}}=\frac{32}{\sqrt{32\cdot 8}}=2$ and $\frac{\pi_1}{\pi_2}=\frac{w_{12}}{w_{21}}=\frac{32}{8}=4$. Thus, the strength ratio corresponding to a given score ratio (as opposed to win ratio) depends on the number of draws. This seems like an undesirable property, though it might be argued that although one cannot directly compare player strength values between different data sets, they do provide meaningful results for a given data set (particularly Consistency, correct ranking for balanced designs). Property~\ref{Luce_Condition}, $\frac{\pi_1}{\pi_2}=\frac{p_{ij}}{p_{ji}}$, satisfied by Davidson's model, has been called an IIA property (Independence of Irrelevant Alternatives)~\cite{Glickman2025}, and indeed, it derives from Luce's choice axiom in the context of comparison of multiple ($>2$) objects~\cite[p.10]{David1988}\cite{Davidson1970}, but when we are considering the effect of draws between the {\em same} set of players as the wins and losses, then the number of draws should, in our view, not be considered {\em irrelevant} alternatives. We return to this assumption of Property~\ref{Luce_Condition} in Section~\ref{sec:reassessment}.

\subsection{Other models}
There are other models for draws, which we do not consider in depth here. For example, draws are included by Rao and Kupper~\cite{Rao1967} by means of a parameter $\theta\ge 1$ such that
$$p_{ij}=\frac{\pi_i}{\pi_i+\theta\pi_j}\,\quad\mbox{and}\quad d_{ij}=(\theta^2-1)p_{ij}p_{ji}=\frac{(\theta^2-1)\pi_i\pi_j}{(\pi_i+\theta\pi_j)(\theta\pi_i+\pi_j)}.$$
Davidson points out that the Rao-Kupper model does not always satisfy Consistency, while his model does (see below).

Glickman's 1999 model~\cite{Glickman1999} deals with dynamic paired comparison experiments, a problem not considered here. We note, however, that for maximum likelihood estimation, it simply treats a draw as equivalent to a win followed by a loss, each counted as half a game, without a draw propensity parameter. We return to this point in Section~\ref{sec:reassessment}.


\section{An alternative model with draws} \label{sec:alternative}
David's assertion that the probability of draws should be accounted for in estimating player strength parameters suggests a natural-looking modification of Davidson's assumption of Property~\ref{Luce_Condition}, and consequently (\ref{bt_davidson}), while accepting Davidson's assumption~(\ref{DavidsonDraw}).

We introduce strength parameters, $\sigma_i$, for $i=1,\ldots,t$, to satisfy
\begin{property}\label{alt_strength_ratio}
\begin{equation}\frac{\sigma_i}{\sigma_j}=\frac{p_{ij}+\frac 12 d_{ij}}{p_{ji}+\frac 12 d_{ij}}\,.\label{strength_ratio}\end{equation}
\end{property}
Property~\ref{alt_strength_ratio} implies that
\begin{equation}\left(\frac{p_{ij}+\frac 12 d_{ij}}{p_{ji}+\frac 12 d_{ij}}\right)\left(\frac{p_{jk}+\frac 12 d_{jk}}{p_{kj}+\frac 12 d_{jk}}\right)=\frac{p_{ik}+\frac 12 d_{ik}}{p_{ki}+\frac 12 d_{ik}}\,.\label{bt_edwards}\end{equation}

This has the advantage over Davidson's model of ensuring that strength parameter ratios depend on expected score ratios, not just the probability of wins (exclusive of draws). On the other hand it has the disadvantage (perhaps surprisingly) that ratings can in some situations not follow the ranking given by scores in balanced designs with three or more players. We demonstrate this fact in Section~\ref{sec:consistency}.

We remark here that whether Davidson's or the Alternative model is more correct could in principle be checked by comparing data to the assumptions (\ref{bt_davidson}) and (\ref{bt_edwards}).

Given~(\ref{strength_ratio}), we calculate:
$$p_{ij}+\frac 12d_{ij}=\frac{\sigma_i}{\sigma_j}\left(p_{ji}+\frac 12 d_{ij}\right)$$
so that
$$\frac 12d_{ij}\left(\frac{\sigma_i}{\sigma_j}-1\right)=p_{ij}-\frac{\sigma_i}{\sigma_j}p_{ji}$$
and
\begin{equation}d_{ij}=\frac{2\left(\sigma_jp_{ij}-\sigma_ip_{ji}\right)}{\sigma_i-\sigma_j}\,,\label{draw_prob}\end{equation}
as long as $\sigma_i\ne\sigma_j$. 
Then, from~(\ref{DavidsonDraw}) and (\ref{draw_prob}),
\begin{equation}\nu \sqrt{p_{ij}p_{ji}}=\frac{2\left(\sigma_jp_{ij}-\sigma_ip_{ji}\right)}{\sigma_i-\sigma_j}\label{draw_base}\end{equation}
$$\frac{\nu}{2}(\sigma_i-\sigma_j)\sqrt{p_{ij}p_{ji}}=\left(\sigma_jp_{ij}-\sigma_ip_{ji}\right)$$
or
$$\frac{\nu}{2}\left(1-\frac{\sigma_j}{\sigma_i}\right)\sqrt{\frac{p_{ji}}{p_{ij}}}=\left(\frac{\sigma_j}{\sigma_i}-\frac{p_{ji}}{p_{ij}}\right)\,.$$
Thus,
$$\frac{p_{ji}}{p_{ij}}+\frac{\nu}{2}\left(1-\frac{\sigma_j}{\sigma_i}\right)\sqrt{\frac{p_{ji}}{p_{ij}}}-\frac{\sigma_j}{\sigma_i}=0$$
from which,
\begin{equation}\sqrt{\frac{p_{ji}}{p_{ij}}}= \frac{\nu}{4}\left(\frac{\sigma_j}{\sigma_i}-1\right)+\sqrt{\frac{\nu^2}{16}\left(1-\frac{\sigma_j}{\sigma_i}\right)^2+\frac{\sigma_j}{\sigma_i}}=\frac{2}{\nu}\phi\left(\frac{\sigma_j}{\sigma_i}\right)\label{pratio}\end{equation}
where we pick the positive root (the other is negative and therefore not relevant), and where, for later use,
\begin{equation}
\phi(x)\equiv \frac{\nu^2}{8}(x-1)+\frac{\nu}{2}\sqrt{\frac{\nu^2}{16}(x-1)^2+x}\,. \label{phi}
\end{equation}
Although $\phi$ depends on parameter $\nu$, writing $\phi(x;\nu)$ everywhere is cumbersome, so we write $\phi(x)$ and just remember the dependence on $\nu$. Note that $\phi(0)=0$, $\phi(1)=\frac{\nu}{2}$, $\phi(\infty)=\infty$, and $\phi(x)>0$ for all $x>0$. Also, if $\nu=0$, then $\phi(x)\equiv 0$, and if $\nu=2$, then $\phi(x)\equiv x$.

Now, using~(\ref{draw_prob}) in~(\ref{ppd}) gives
$$p_{ij}+p_{ji}+\frac{2(\sigma_jp_{ij}-\sigma_ip_{ji})}{\sigma_i-\sigma_j}=1$$
so
$$(\sigma_i-\sigma_j)(p_{ij}+p_{ji})+2(\sigma_jp_{ij}-\sigma_ip_{ji})=\sigma_i-\sigma_j$$
$$(\sigma_i+\sigma_j)p_{ij}-(\sigma_i+\sigma_j)p_{ji}=\sigma_i-\sigma_j$$
and
\begin{equation}
p_{ji}=p_{ij}+\frac{\sigma_j-\sigma_i}{\sigma_j+\sigma_i}=p_{ij}+\frac{\left(\frac{\sigma_j}{\sigma_i}-1\right)}{\left(\frac{\sigma_j}{\sigma_i}+1\right)}\,.\label{psigma}
\end{equation}
Squaring~(\ref{pratio}) gives
\begin{eqnarray*}
p_{ji} &=& p_{ij}\left(\frac{2}{\nu}\phi\left( \frac{\sigma_j}{\sigma_i} \right)\right)^2=p_{ij}\left(\frac{\nu}{4}\left(\frac{\sigma_j}{\sigma_i}-1\right)+\sqrt{\frac{\nu^2}{16}\left(1-\frac{\sigma_j}{\sigma_i}\right)^2+\frac{\sigma_j}{\sigma_i}}\right)^2 \\
&=& p_{ij}\left(\frac{\nu^2}{8}\left(\frac{\sigma_j}{\sigma_i}-1\right)^2 +\frac{\sigma_j}{\sigma_i}+\frac{\nu}{2}\left(\frac{\sigma_j}{\sigma_i}-1\right)\sqrt{\frac{\nu^2}{16}\left(1-\frac{\sigma_j}{\sigma_i}\right)^2+\frac{\sigma_j}{\sigma_i}} \right) \\
&=& p_{ij}\left[\left(\frac{\sigma_j}{\sigma_i}-1\right)\left(\frac{\nu^2}{8}\left(\frac{\sigma_j}{\sigma_i}-1\right) +\frac{\nu}{2}\sqrt{\frac{\nu^2}{16}\left(1-\frac{\sigma_j}{\sigma_i}\right)^2+\frac{\sigma_j}{\sigma_i}} \right)+\frac{\sigma_j}{\sigma_i}\right] \\
&=& p_{ij}\left[\left(\frac{\sigma_j}{\sigma_i}-1\right)\phi\left(\frac{\sigma_j}{\sigma_i}\right)+\frac{\sigma_j}{\sigma_i}\right]
\label{pratio_squared}
\end{eqnarray*}
and combining with~(\ref{psigma}) gives
$$ p_{ij}\left[\left(\frac{\sigma_j}{\sigma_i}-1\right)\phi\left(\frac{\sigma_j}{\sigma_i}\right)+\frac{\sigma_j}{\sigma_i}-1 \right]= \frac{\left(\frac{\sigma_j}{\sigma_i}-1\right)}{\left(\frac{\sigma_j}{\sigma_i}+1\right)}$$
or
\begin{equation*}
\frac{1}{p_{ij}} =  \left(\frac{\sigma_j}{\sigma_i}+1\right)\left[1+\phi\left(\frac{\sigma_j}{\sigma_i}\right)\right] 
\end{equation*}
or
\begin{equation}
p_{ij}=\left(\frac{\sigma_i}{\sigma_i+\sigma_j}\right)\left(\frac{1}{1+\phi\left(\frac{\sigma_j}{\sigma_i}\right)}\right)
=\left(\frac{\sigma_i}{\sigma_i+\sigma_j}\right)\left(\frac{\sigma_i}{\sigma_i+\sigma_i\phi\left(\frac{\sigma_j}{\sigma_i}\right)}\right)
\,,\label{p_by_sigmaratio}
\end{equation}
which expresses $p_{ij}$ solely in terms of strength parameters $\sigma_i$ and $\sigma_j$ (as well as parameter $\nu$).

If $\sigma_i=\sigma_j$, then $p_{ij}=p_{ji}$ and $d_{ij}=1-2p_{ij}$, so $\nu p_{ij}=1-2p_{ij}$ and $p_{ij}=\frac{1}{2+\nu}$, which agrees with Equation~\eqref{p_by_sigmaratio}.


Also, directly from~(\ref{strength_ratio}) and~(\ref{DavidsonDraw}),
\begin{equation} \frac{\sigma_i}{\sigma_j} = \frac{p_{ij}+\frac{\nu}{2}\sqrt{p_{ij}p_{ji}}}{p_{ji}+\frac{\nu}{2}\sqrt{p_{ij}p_{ji}}} =\frac{ \sqrt{p_{ij}}(\sqrt{p_{ij}}+\frac{\nu}{2}\sqrt{p_{ji}}) }{ \sqrt{p_{ji}}(\sqrt{p_{ji}}+\frac{\nu}{2}\sqrt{p_{ij}})   } = \frac{ \sqrt{\frac{p_{ij}}{p_{ji}}} +\frac{\nu}{2} }{ \sqrt{\frac{p_{ji}}{p_{ij}}} +\frac{\nu}{2} }\,.\label{SigmaRatio2}\end{equation}

\begin{figure}
    \centering
    \includegraphics[scale=.5]
    {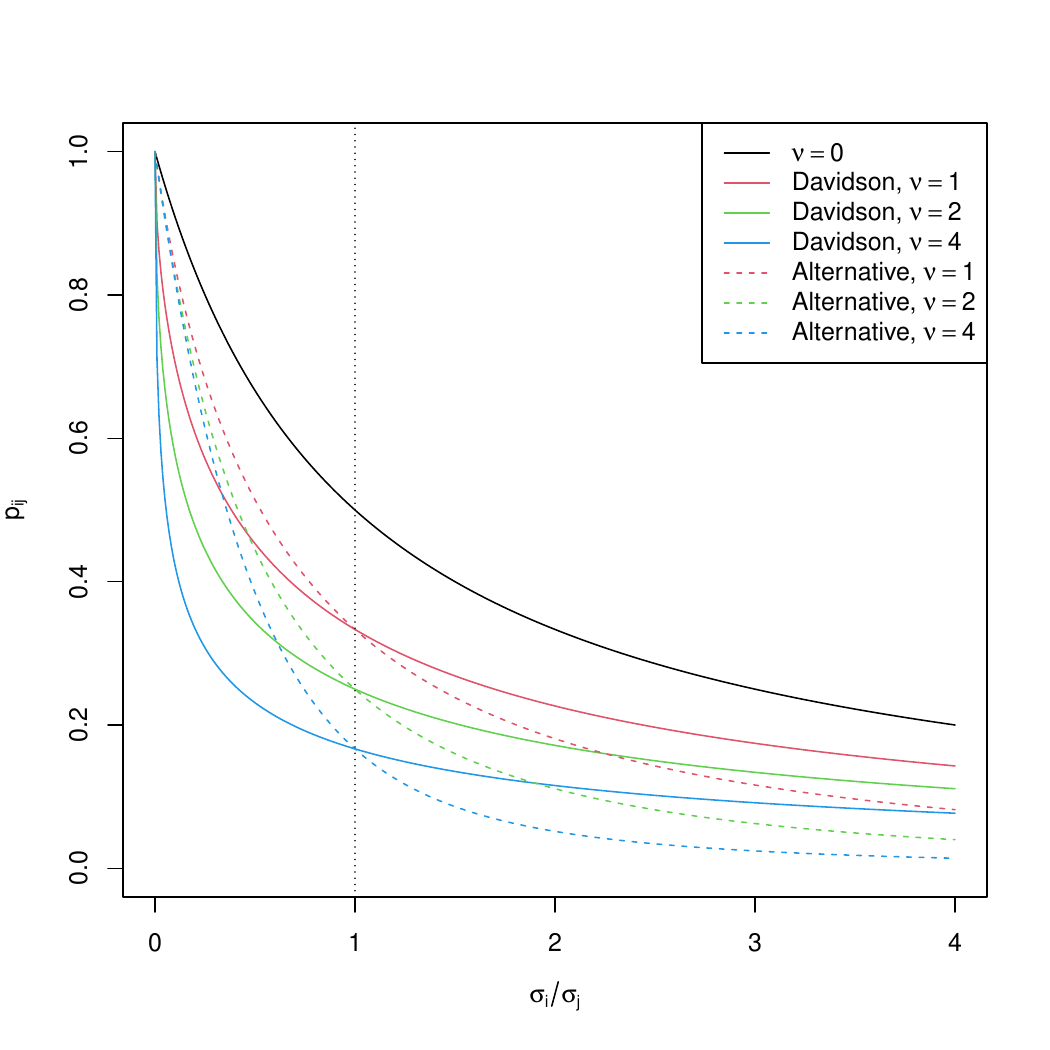}
    \caption{Win probabilities in the Davidson and Alternative model with various values of the draw-propensity parameter $\nu$.}
    \label{fig:tied_pools_plot}
\end{figure}

Figure~\ref{fig:tied_pools_plot} shows the form of $p_{ij}$ for a range of values of $\frac{\sigma_1}{\sigma_j}$ and various values of $\nu$ in both the Davidson model and the Alternative model.

\subsection{Special cases of the Alternative model: $\nu=0$ and $\nu=2$}
Note that when $\nu=0$, the Alternative model reduces to the Bradley-Terry model. 

Another special case is important. When
$\nu=2$, Equation~\eqref{SigmaRatio2} becomes
\begin{equation}\frac{\sigma_i}{\sigma_j}=\sqrt{\frac{p_{ij}}{p_{ji}}}\,.\label{SigmaRatioSimple}\end{equation}
Thus, the ratio of strength parameters is the square root of the ratio for the Davidson model, {\em i.e.}, 
\begin{equation}\frac{\sigma_i}{\sigma_j}=\sqrt{\frac{\pi_i}{\pi_j}}\,. \label{sigma_pi}
\end{equation}
Also, since $\phi(x)=x$, Equation~\eqref{p_by_sigmaratio} becomes
\begin{equation}p_{ij} = \left( \frac{\sigma_i}{\sigma_i+\sigma_j} \right)^2\,.\label{pijEdwardsSimple}\end{equation}
Thus, the probability of a win looks like the probability of two wins in the Bradley-Terry model without draws. In effect, this model with $\nu=2$ is what one would obtain from the Bradley-Terry model without draws if each game were considered to be {\em two} comparisons instead of one, with a draw corresponding to a result of 10 or 01, and a win (loss) corresponding to a result of 11 (or 00). This is a natural way of capturing the fact that a draw is considered equivalent (in scoring) to half a win, {\em i.e.}, two draws are equivalent to a win and a loss. This also implies that variances are halved in relation to the Bradley-Terry model, reflecting the finer information available when draws are possible.


In fact, the Davidson model with $\nu=2$ also gives
\begin{equation} p_{ij}=\frac{\pi_i}{\pi_i+\pi_j+2\sqrt{\pi_i\pi_j}} = \left(\frac{\sqrt{\pi}_i}{\sqrt{\pi_i}+\sqrt{\pi_j}}\right)^2\label{pijDavidsonSimple}\end{equation}
which is equivalent to the above if $\sigma_i=\sqrt{\pi_i}$ and $\sigma_j=\sqrt{\pi_j}$, and is consistent with the observation in Equation~\eqref{sigma_pi} above, since
$$p_{ij}=\left(\frac{1}{1+\frac{\sigma_j}{\sigma_i}}\right)^2=\left(\frac{1}{1+\sqrt{\frac{\pi_j}{\pi_i}}}\right)^2\,.$$

\section{Maximum likelihood estimation}\label{sec:max_likelihood}

\subsection{Davidson model}
The likelihood of a data set in the Davidson model is
\begin{eqnarray*}
L &=& \prod_{i,j} p_{ij}^{w_{ij}} \prod_{i<j}d_{ij}^{t_{ij}} = \prod_{i,j} \left( \frac{\pi_i}{\pi_i+\pi_j+\nu\sqrt{\pi_i \pi_j}} \right)^{w_{ij}} \times \prod_{i<j} \left( \frac{\nu\sqrt{\pi_i \pi_j} }{\pi_i+\pi_j+\nu\sqrt{\pi_i \pi_j}} \right)^{t_{ij}} \\
&=& \frac{\nu^T \prod_{i<j}\pi_i^{w_{ij}+\frac{1}{2} t_{ij}} \prod_{i>j} \pi_i^{w_{ij}+\frac{1}{2}t_{ij}}}{ \prod_{i<j} (\pi_i+\pi_j+\nu\sqrt{\pi_i \pi_j})^{ n_{ij} } } = \frac{\nu^T\prod_{i,j}\pi_i^{s_{ij} }} { \prod_{i<j} (\pi_i+\pi_j+\nu\sqrt{\pi_i \pi_j})^{n_{ij}} } =  \frac{\nu^T\prod_{i}\pi_i^{s_{i} }} { \prod_{i<j}(\pi_i+\pi_j+\nu\sqrt{\pi_i \pi_j})^{n_{ij}} }
\end{eqnarray*}
and the log likelihood is
$$ \log(L) = T\log(\nu) + \sum_i s_i\log(\pi_i) - \sum_{i<j} n_{ij} \log(\pi_i+\pi_j+\nu\sqrt{\pi_i \pi_j})\,. $$

To find the maximum of the log-likelihood function, we differentiate with respect to each parameter and set these to $0$.
\begin{eqnarray*}
 \frac{\partial}{\partial \pi_k}\log(L) &=& \frac{s_k}{\pi_k} - \sum_{i<k} \frac{n_{ik}\left(1+\frac{\nu}{2}\sqrt{\frac{\pi_i}{\pi_k}}\right)}{\pi_i+\pi_k+\nu\sqrt{\pi_i\pi_k}} + \sum_{j>k} \frac{n_{kj}\left(1+\frac{\nu}{2}\sqrt{\frac{\pi_j}{\pi_k}}\right)}{\pi_k+\pi_j+\nu\sqrt{\pi_k\pi_j}}  \\
 &=&  \frac{s_k}{\pi_k} - \sum_{i} \frac{n_{ik}\left(1+\frac{\nu}{2}\sqrt{\frac{\pi_i}{\pi_k}}\right)}{\pi_i+\pi_k+\nu\sqrt{\pi_i\pi_k}} =\frac{s_k}{\pi_k} - G_k(\bm{\pi},\nu) =0
 \end{eqnarray*}
letting $\bm{\pi}=(\pi_1,\ldots,\pi_t)^\top$, and where we use the fact that $n_{kj}=n_{jk}$ (and of course $n_{kk}=0$). And,
\begin{eqnarray*}
  \frac{\partial}{\partial \nu}\log(L)&=& \frac{T}{\nu} - \sum_{i<j}\frac{n_{ij}\sqrt{\pi_i\pi_j}}{\pi_i+\pi_j+\nu\sqrt{\pi_i\pi_j}} \\
  &=&  \frac{T}{\nu} - \frac 12 \sum_{i,j}\frac{n_{ij}\sqrt{\pi_i\pi_j}}{\pi_i+\pi_j+\nu\sqrt{\pi_i\pi_j}} =\frac{T}{\nu}-H(\bm{\pi},\nu) = 0\,.
  \end{eqnarray*}
The maximum likelihood solution is found by iteratively solving the system
\begin{subequations}
\label{iteration_Davidson}
\begin{align}
  \pi_k &\leftarrow \frac{s_k}{G_k(\bm{\pi},\nu)}\,,\quad k=1,\ldots,t, \\
  \nu &\leftarrow \frac{T}{H(\bm{\pi},\nu)  } \,,
\end{align}
\end{subequations}
along with a normalization criterion, like $\sum_k \pi_k = 1$. 
In the above system of $t+1$ equations, the number of parameters is $t+1$, since there are $t$ strength parameters $\pi_k, k=1,\ldots, t$, as well as $\nu$, but since the strength parameter vector, $\bm{\pi}=(\pi_1,\ldots,\pi_t)^\top$, can be scaled arbitrarily without changing scores or numbers of draws, we must lose one degree of freedom by such a normalization. Then the number of independent parameters is $t$. Indeed, one could rewrite the system in terms of parameter ratios, since $s_k=\frac{\pi_k}{\pi_t}G_k\left(\frac{\bm{\pi}}{\pi_t},\nu\right)$, which makes it clear that there are only $t$ independent parameters, the $t-1$ ratios, $\frac{\pi_k}{\pi_t}$ for $k=1,\ldots,t-1$ (but not $\frac{\pi_t}{\pi_t}$, which of course is fixed at $1$). Now it seems there are $t+1$ equations in $t$ unknowns, but one of the first $t$ equations is dependent on the others. It is clear that if all data ($s_k$, $n_{ik}$, and $T$), are scaled uniformly by a factor $\lambda>0$, the equations are the same, so the parameters estimated must be the same. In effect, only the score ratios are needed, not the scores themselves. So the first $t-1$ equations could be written
  $$ \frac{s_k}{s_t}=\frac{\pi_k}{\pi_t}\frac{G_t\left(\frac{\bm{\pi}}{\pi_t},\nu\right)}{G_k\left(\frac{\bm{\pi}}{\pi_t},\nu\right)}\,,\quad k=1,\ldots, t-1\,. $$
In effect, there are $t$ independent equations (for the $s$-ratios and $T$) in $t$ independent unknowns (the $\pi$-ratios and $\nu$).
  
  Also, note that in a balanced tournament, where all $n_{ij}\equiv n$, the number of data points in the set of equations for the maximum log-likelihood is $t+2$, since they are the components of the vector $s=(s_1,\ldots,s_t)$ as well as $n$ and $T$. But $n$ is not independent of $s$, since $\sum_k s_k = \frac{n}{2}t(t-1)$, and an arbitrary scaling of all $s_{ij}=w_{ij}+\frac 12 t_{ij}$ scales $s$ and $T$ (and $n$) by the same factor, $\lambda$, and does not change the equations, so we lose another degree of freedom by fixing this scaling, and the number of independent data points is $t$. A convenient normalization here is to make $n=1$, so that each $s_{ij}+s_{ji}=1$, and $\sum_k s_k = \frac 12 t(t-1)$. This implies that the reverse problem, to find the data values given the parameters, is also determined ($t$ independent equations in $t$ unknowns) for balanced tournaments. This is not the case, of course, for unbalanced tournaments, since then there are many more independent data values, $n_{ik}$, so many sets of data values will be consistent with the same set of parameters.
  
  The Consistency property, the fact that in a balanced tournament (where $n_{ij}=n$, for all $i\ne j$), strength parameters respect the order of scores, is shown by Davidson~\cite{Davidson1970} as follows.
  \begin{prop} The Davidson model is consistent.\end{prop}
  {\bf Proof:}
  \begin{eqnarray*}
  s_j-s_k &=& \pi_jG_j(\bm{\pi},\nu)-\pi_kG_k(\bm{\pi},\nu) = \pi_j n \sum_{i\ne j} \frac{\left(1+\frac{\nu}{2}\sqrt{\frac{\pi_i}{\pi_j}}\right)}{\pi_i+\pi_j+\nu\sqrt{\pi_i\pi_j}} - \pi_k n \sum_{i\ne k} \frac{\left(1+\frac{\nu}{2}\sqrt{\frac{\pi_i}{\pi_k}}\right)}{\pi_i+\pi_k+\nu\sqrt{\pi_i\pi_k}} \\
  &=& n \sum_{i\ne j} \frac{\left(\pi_j+\frac{\nu}{2}\sqrt{\pi_i\pi_j}\right)}{\pi_i+\pi_j+\nu\sqrt{\pi_i\pi_j}} - n\sum_{i\ne k} \frac{\left(\pi_k+\frac{\nu}{2}\sqrt{\pi_i \pi_k}\right)}{\pi_i+\pi_k+\nu\sqrt{\pi_i\pi_k}} \\
  &=& n \sum_{i} \frac{\left(\pi_j+\frac{\nu}{2}\sqrt{\pi_i\pi_j}\right)}{\pi_i+\pi_j+\nu\sqrt{\pi_i\pi_j}} - n\sum_{i} \frac{\left(\pi_k+\frac{\nu}{2}\sqrt{\pi_i \pi_k}\right)}{\pi_i+\pi_k+\nu\sqrt{\pi_i\pi_k}}- n\frac{\pi_j+\frac{\nu}{2}\sqrt{\pi_j\pi_j}}{\pi_j+\pi_j+\nu\sqrt{\pi_j\pi_j}}+ n\frac{\pi_k+\frac{\nu}{2}\sqrt{\pi_k\pi_k}}{\pi_k+\pi_k+\nu\sqrt{\pi_k\pi_k}} \\
  &=& n \sum_{i} \frac{\left(\pi_j+\frac{\nu}{2}\sqrt{\pi_i\pi_j}\right)\left(\pi_i+\pi_k+\nu\sqrt{\pi_i\pi_k}\right) - \left(\pi_k+\frac{\nu}{2}\sqrt{\pi_i \pi_k}\right)\left(\pi_i+\pi_j+\nu\sqrt{\pi_i\pi_j}\right)}{\left(\pi_i+\pi_j+\nu\sqrt{\pi_i\pi_j}\right)\left(\pi_i+\pi_k+\nu\sqrt{\pi_i\pi_k}\right) } -\frac{n}{2} +\frac{n}{2}\\
  &=& n \sum_{i} \frac{\left(\pi_i(\pi_j-\pi_k)+\frac{\nu}{2}\pi_j\sqrt{\pi_i \pi_k}-\frac{\nu}{2}\pi_k\sqrt{\pi_i \pi_j}+\frac{\nu}{2}\pi_i\sqrt{\pi_i \pi_j}-\frac{\nu}{2}\pi_i\sqrt{\pi_i \pi_k}\right)}{\left(\pi_i+\pi_j+\nu\sqrt{\pi_i\pi_j}\right)\left(\pi_i+\pi_k+\nu\sqrt{\pi_i\pi_k}\right) } \\
 &=& n \sum_{i} \frac{\left(\pi_i(\sqrt{\pi_j}-\sqrt{\pi_k})(\sqrt{\pi_j}+\sqrt{\pi_k})+(\sqrt{\pi_j}-\sqrt{\pi_k})\frac{\nu}{2}\sqrt{\pi_i \pi_j \pi_k}+(\sqrt{\pi_j}-\sqrt{\pi_k})\frac{\nu}{2}\pi_i^{3/2}\right)}{\left(\pi_i+\pi_j+\nu\sqrt{\pi_i\pi_j}\right)\left(\pi_i+\pi_k+\nu\sqrt{\pi_i\pi_k}\right) } \\
  &=& n(\sqrt{\pi_j}-\sqrt{\pi_k}) \sum_{i} \frac{\left(\pi_i(\sqrt{\pi_j}+\sqrt{\pi_k})+\frac{\nu}{2}\sqrt{\pi_i \pi_j \pi_k}+\frac{\nu}{2}\pi_i^{3/2}\right)}{\left(\pi_i+\pi_j+\nu\sqrt{\pi_i\pi_j}\right)\left(\pi_i+\pi_k+\nu\sqrt{\pi_i\pi_k}\right) } \\
 \end{eqnarray*}
 and since the sum in the last expression is positive, $(s_j-s_k)$ and $(\sqrt{\pi_j}-\sqrt{\pi_k})$ have the same sign, and so, therefore, too, does $(\pi_j-\pi_k)$. \qed
 
 Note that the same proof works (more simply) for the Bradley-Terry model (without draws), by setting $\nu=0$.
 
 Note also that in a 2-player tournament, where $t=2$, we have $n_{12}=n_{21}=n$ and
 $$ s_1=\pi_1 G_1(\pi,\nu)=\frac{n\left(1+\frac{\nu}{2}\sqrt{\frac{\pi_2}{\pi_1}}\right)}{\pi_1+\pi_2+\nu\sqrt{\pi_1\pi_2}} $$
and similarly for $s_2$, so
$$ \frac{s_1}{s_2}=\frac{\pi_1 G_1(\pi,\nu)}{\pi_2 G_2(\pi,\nu)}=\frac{\pi_1+\frac{\nu}{2}\sqrt{\pi_1\pi_2}}{\pi_2+\frac{\nu}{2}\sqrt{\pi_1\pi_2}} =\frac{p_{12}+\frac 12 d_{12}}{p_{21}+\frac 12 d_{21}} \, ,$$
and, of course, by definition, 
$$ \frac{s_1}{s_2}=\frac{w_{12}+\frac 12 t_{12}}{w_{21}+\frac 12 t_{21}}\,.$$

\subsection{Alternative model} \label{sec:alternative_likelihood}

To find the maximum likelihood estimate for the Alternative model, we need $p_{ij}$ expressed in terms of $\sigma_i$, $\sigma_j$ and $\nu$, which is given by Equation~\eqref{p_by_sigmaratio}.

We have 
$$L= \prod_{i,j} p_{ij}^{w_{ij}} \prod_{i<j}d_{ij}^{t_{ij}} = \left(\prod_{i,j}p_{ij}^{w_{ij}}\right)\nu^T\left(\prod_{i<j}\sqrt{p_{ij}}^{t_{ij}}\right)\left(\prod_{i<j}\sqrt{p_{ji}}^{t_{ji}}\right) = \nu^T\prod_{i,j}p_{ij}^{s_{ij}}$$
so
\begin{eqnarray*}
 \log(L) &=& T\log(\nu)+\sum_{i,j}s_{ij}\log(p_{ij}) \\
 &=& T\log(\nu)+\sum_{i,j}s_{ij} \left[ \log(\sigma_i)-\log(\sigma_i+\sigma_j)-\log\left(1+\phi\left(\frac{\sigma_j}{\sigma_i}\right)\right) \right]  \\
 &=& T\log(\nu)+\sum_{i}s_{i} \log(\sigma_i) -\sum_{i,j}s_{ij} \log(\sigma_i+\sigma_j)-\sum_{i,j}s_{ij}\log\left(1+\phi\left(\frac{\sigma_j}{\sigma_i}\right)\right)    \,.
 \end{eqnarray*}

 Then the maximum likelihood is obtained from 
\begin{subequations}
\label{max_likelihood_alternative}
\begin{align}
0 = \frac{\partial}{\partial \sigma_k}\log(L) &= \frac{s_k}{\sigma_k} - \sum_j \frac{s_{kj}}{\sigma_k+\sigma_j} - \sum_i \frac{s_{ik}}{\sigma_i+\sigma_k} 
- \sum_j \frac{s_{kj}\left(-\frac{\sigma_j}{\sigma_k^2}\phi'\left(\frac{\sigma_j}{\sigma_k}\right)\right)}{1+\phi\left(\frac{\sigma_j}{\sigma_k}\right)} - \sum_i \frac{s_{ik}\frac{1}{\sigma_i}\phi'\left(\frac{\sigma_k}{\sigma_i}\right)}{1+\phi\left(\frac{\sigma_k}{\sigma_i}\right)} \notag \\
&=  \frac{s_k}{\sigma_k} - \sum_i \frac{n_{ik}}{\sigma_i+\sigma_k} + \sum_i \frac{s_{ki}\left(\frac{\sigma_i}{\sigma_k^2}\phi'\left(\frac{\sigma_i}{\sigma_k}\right)\right)}{1+\phi\left(\frac{\sigma_i}{\sigma_k}\right)} - \sum_i \frac{s_{ik}\frac{1}{\sigma_i}\phi'\left(\frac{\sigma_k}{\sigma_i}\right)}{1+\phi\left(\frac{\sigma_k}{\sigma_i}\right)} =\frac{s_k}{\sigma_k}-G_k(\bm{\sigma},\nu) \label{sG}\\
0 = \frac{\partial}{\partial\nu}\log(L) &= \frac{T}{\nu} - \sum_{i,j} \frac{s_{ij}\phi_{\nu}\left(\frac{\sigma_j}{\sigma_i}\right)}{1+\phi\left(\frac{\sigma_j}{\sigma_i}\right)} =\frac{T}{\nu}-H(\bm{\sigma},\nu) \,, \label{TH}
\end{align}
\end{subequations}
where
\begin{eqnarray*} G_k(\bm{\sigma},\nu) &=& \sum_i \frac{n_{ik}}{\sigma_i+\sigma_k} - \sum_i \frac{s_{ki}\left(\frac{\sigma_i}{\sigma_k^2}\phi'\left(\frac{\sigma_i}{\sigma_k}\right)\right)}{1+\phi\left(\frac{\sigma_i}{\sigma_k}\right)} + \sum_i \frac{s_{ik}\frac{1}{\sigma_i}\phi'\left(\frac{\sigma_k}{\sigma_i}\right)}{1+\phi\left(\frac{\sigma_k}{\sigma_i}\right)}  \\
H(\bm{\sigma},\nu) &=& \sum_{i,j} \frac{s_{ij}\phi_{\nu}\left(\frac{\sigma_j}{\sigma_i}\right)}{1+\phi\left(\frac{\sigma_j}{\sigma_i}\right)}
\end{eqnarray*}
which suggests the iteration
\begin{subequations}
\label{iteration_alternative}
\begin{align}
    \sigma_k &\leftarrow \frac{s_k}{G_k(\bm{\sigma},\nu)}\,,\quad k=1,\ldots,t, \\
    \nu &\leftarrow \frac{T}{H(\bm{\sigma},\nu)}\,.
    \end{align}
\end{subequations}

Here, 
\begin{equation}\phi'(x)=\frac{\nu^2}{8}+\frac{\frac{\nu}{4}\left[\frac{\nu^2}{8}(x-1)+1\right]}{\sqrt{\frac{\nu^2}{16}(x-1)^2+x}}=\frac{\nu}{4q(x)}\left[\frac{\nu q(x)}{2}+\frac{\nu^2}{8}(x-1)+1\right]\,, \label{phip}
\end{equation}
where 
\begin{equation} 
q(x) \equiv \sqrt{\frac{\nu^2}{16}(x-1)^2+x}\,. \label{q}
\end{equation}
Also,
\begin{eqnarray*}
\phi_{\nu}(x) &=& \frac{\nu}{4}(x-1)+\frac{1}{2}\sqrt{\frac{\nu^2}{16}(x-1)^2+x}+\frac{\frac{\nu}{4}\left[\frac{\nu}{8}(x-1)^2\right]}{\sqrt{\frac{\nu^2}{16}(x-1)^2+x}} \\
&=&  \frac{\nu}{4}(x-1) + \frac{ \frac{1}{2}\left[\frac{\nu^2}{16}(x-1)^2+x\right] +\frac{\nu^2}{32}(x-1)^2 }{\sqrt{\frac{\nu^2}{16}(x-1)^2+x}} \\
&=&  \frac{\nu}{4}(x-1) + \frac{ \frac{\nu^2}{16}(x-1)^2+\frac{x}{2}  }{\sqrt{\frac{\nu^2}{16}(x-1)^2+x}} \,,
\end{eqnarray*}
so 
$$\phi'\left(\frac{\sigma_j}{\sigma_i}\right)=\frac{\nu^2}{8}+\frac{\frac{\nu}{4}\left[\frac{\nu^2}{8}\left(\frac{\sigma_j}{\sigma_i}-1\right)+1\right]}{\sqrt{\frac{\nu^2}{16}\left(\frac{\sigma_j}{\sigma_i}-1\right)^2+\frac{\sigma_j}{\sigma_i}}}\,,\quad 
\phi_{\nu}\left(\frac{\sigma_j}{\sigma_i}\right) =  \frac{\nu}{4}\left(\frac{\sigma_j}{\sigma_i}-1\right) + \frac{ \frac{\nu^2}{16}\left(\frac{\sigma_j}{\sigma_i}-1\right)^2+\frac{\sigma_j}{2\sigma_i} }{\sqrt{\frac{\nu^2}{16}\left(\frac{\sigma_j}{\sigma_i}-1\right)^2+\frac{\sigma_j}{\sigma_i}}}  \,.$$
Note that $\phi'(0)=1$, $\phi'(1)=\frac{\nu^2}{8}+\frac{\nu}{4}$, $\phi'(\infty)=\frac{\nu^2}{4}$, and if $\nu=2$ then $\phi'(x)\equiv 1$. It can also be shown that $\phi'(x)>0$ (except when $\nu=0$), so $\phi(x)$ is increasing, as follows. 
Since $\phi'(0)=1$, we need only show that there is no $x$ for which $\phi'(x)=0$. Suppose there is. From Equation~\eqref{phip}, $\phi'(x)=0$ when 
$$\frac{\nu}{2}q(x)=\frac{\nu^2}{8}(1-x)-1\,,$$
noting that $q(x)>0$ for all $x\ge 0$.
Squaring both sides, we find that if $\phi'(x)=0$ then
$$ \frac{\nu^2}{4}\left[\frac{\nu^2}{16}(x-1)^2+x\right]=\frac{\nu^2}{64}(1-x)^2-\frac{\nu^2}{4}(1-x)+1\,,$$
which reduces to $0=1-\frac{\nu^2}{4}$, or $\nu=2$. However, when $\nu=2$, we have $\phi'(x)\equiv 1$, so we have a contradiction, and there is no $x$ for which $\phi'(x)=0$. 


In practice, convergence of the iteration is improved if instead of using Equation~\eqref{sG}, we write 
$$ 0=\frac{3s_k}{\sigma_k}-\left(G_k(\bm{\sigma},\nu)+\frac{2s_k}{\sigma_k}\right)$$
and iterate using
\begin{subequations}
\label{iteration_alternative_practical}
\begin{align}
\sigma_k &\leftarrow \frac{3s_k}{\left(G_k(\bm{\sigma},\nu)+\frac{2s_k}{\sigma_k}\right)} \\
\nu &\leftarrow \frac{T}{H(\bm{\sigma},\nu)}\,.
\end{align}
\end{subequations}
If this iteration is not guaranteed to converge, one could also use Newton's method, and we do not pursue this question further.

\section{Consistency in the Alternative model}\label{sec:consistency}
We now wish to demonstrate that in balanced designs, Property~\ref{consistency_property} does not necessarily hold, that is, strength parameters {\bf do not necessarily} respect the order of scores, if there are $t>2$ players. 

Whether or not the iteration converges, we can attempt to solve the reverse problem for $s_{ij}$ and $T$ given values for parameters $(\bm{\sigma},\nu)$, where $\bm{\sigma}=(\sigma_1,\ldots,\sigma_t)$, since this is then just a linear system. However, when there are $t\ge 3$ players, this is now an underdetermined system, so there are many possible solutions. For a $t$-player system, there are $t+1$ parameters in $(\bm{\sigma},\nu)$, but we must subtract 1 because only the ratios of player strengths are determined ({\em i.e.}, $\bm{\sigma}$ can be scaled arbitrarily without changing scores), so there are only $t$ independent parameters after a normalization (like $\sum_i \sigma_i=1$ or $\sigma_t=1$). 

In the Davidson model, the individual $s_{ij}$ do not appear in the set of equations for the maximum log-likelihood, only the $t$ values of $s_i$. Thus, for a balanced tournament, the set of data values to solve for in the inverse problem, consists of $s$, $n$, and $T$, a total of $t+2$ values. As discussed above, however, there are two additional conditions, from the fact that $\sum_i s_i = \sum_{i<j} n_{ij} = \frac{n}{2}t(t-1)=N$, and from the fact that all scores $s_{ij}=w_{ij}+\frac 12 t_{ij}$ can be scaled arbitrarily for any set of parameters, so a normalization can be imposed, such as $N=\frac 12 t(t-1)$, which for a balanced tournament implies $n=1$. Thus, in the Davidson model, there are only $t$ independent data values to solve for, 
and there are $t$ independent equations ($t+1$ equations but one of the equations for $s_i$ must be dependent on the others, since only ratios $\frac{s_i}{s_j}$ are determined). 

In the Alternative model, the pairwise scores, $s_{ij}$, do appear in the set of equations for the maximum log-likelihood, so the set of data values to solve for in the inverse problem for a balanced tournament comprises the $\frac 12 t(t-1)$ values $s_{ij}$ for $i<j$, as well as $n$ and $T$, in total $\frac 12 t(t-1)+2$ of them. Here,  $n$ is not independent of the score parameters, $s_{ij}$ for $i<j$. Indeed, one needs $n$ and these scores to determine each $s_{ji}$ for $i<j$. There is still one additional constraint, since an arbitrary uniform scaling can be applied to all scores (thus, scaling $n$ and $T$ by the same factor), so we can normalize by choosing $n=1$, leaving us with $\frac 12 t(t-1)+1$ independent data values to determine.

When $t=2$, the linear system is exactly determined, since there are $t=2$ independent parameters and $\frac 12 t(t-1)+1=2$ independent data values, and when $t=3$, it is already underdetermined, since there are $t=3$ independent parameters and $\frac 12 t(t-1)+1=4$ data values. We can show inconsistency by showing that there exist parameter ratios $\frac{\sigma_i}{\sigma_j}$ and score ratios $\frac{s_i}{s_j}$ such that their logarithms have opposite sign (one ratio is bigger than $1$ and the other ratio is smaller than $1$).

\subsection{Two-player tournaments}
If $t=2$, then $s_{12}+s_{21}=n$, and as mentioned above, we can fix $n=1$. Also $s_1=s_{12}$ and $s_2=s_{21}=1-s_1$. And
\begin{eqnarray*}
s_1 &=& \sigma_1 G_1(\sigma,\nu)=\sigma_1\left[\frac{1}{\sigma_1+\sigma_2}-\frac{s_{12}\left(\frac{\sigma_2}{\sigma_1^2}\phi'\left(\frac{\sigma_2}{\sigma_1}\right)\right)}{1+\phi\left(\frac{\sigma_2}{\sigma_1}\right)} + \frac{s_{21}\frac{1}{\sigma_2}\phi'\left(\frac{\sigma_1}{\sigma_2}\right)}{1+\phi\left(\frac{\sigma_1}{\sigma_2}\right)}\right] \\
&=& \frac{1}{1+\frac{\sigma_2}{\sigma_1}} - s_1\frac{ \frac{\sigma_2}{\sigma_1}\phi'\left(\frac{\sigma_2}{\sigma_1}\right)}{1+\phi\left(\frac{\sigma_2}{\sigma_1}\right)} + (1-s_1)\frac{\frac{\sigma_1}{\sigma_2}\phi'\left(\frac{\sigma_1}{\sigma_2}\right)}{1+\phi\left(\frac{\sigma_1}{\sigma_2}\right)} \,,
\end{eqnarray*}
\begin{equation}
s_1 \left[1+\frac{ \frac{\sigma_2}{\sigma_1}\phi'\left(\frac{\sigma_2}{\sigma_1}\right)}{1+\phi\left(\frac{\sigma_2}{\sigma_1}\right)} +\frac{\frac{\sigma_1}{\sigma_2}\phi'\left(\frac{\sigma_1}{\sigma_2}\right)}{1+\phi\left(\frac{\sigma_1}{\sigma_2}\right)} \right]  =  \frac{1}{1+\frac{\sigma_2}{\sigma_1}}  + \frac{\frac{\sigma_1}{\sigma_2}\phi'\left(\frac{\sigma_1}{\sigma_2}\right)}{1+\phi\left(\frac{\sigma_1}{\sigma_2}\right)} \,.
\label{eq:sss}
\end{equation}
Note that when $\frac{\sigma_1}{\sigma_2}=1$, we have $\phi(1)=\frac{\nu}{2}$ and $\phi'(1)=\frac{\nu^2}{8}+\frac{\nu}{4}=\frac{\nu}{4}\left(\frac{\nu}{2}+1\right)$, so 
$$s_1\left[1+2\frac{\frac{\nu}{4}\left(\frac{\nu}{2}+1\right)}{1+\frac{\nu}{2}}\right] = \frac{1}{2}+\frac{\frac{\nu}{4}\left(\frac{\nu}{2}+1\right)}{1+\frac{\nu}{2}} $$
$$s_1 = \frac{\frac 12 +\frac{\nu}{4}}{1+\frac{\nu}{2}} =\frac 12\,,\quad \mbox{so } s_2=\frac 12\,. $$

\begin{prop} \label{prop:2player}
In a (necessarily balanced) 2-player tournament, that is when $t=2$, we have
$\frac{s_1}{s_2}=\frac{\sigma_1}{\sigma_2}$.
\end{prop}

{\bf Proof}:
From Equation~\eqref{eq:sss},
\begin{eqnarray*} \frac{s_1}{s_2} &=& \frac{  \frac{1}{1+\frac{\sigma_2}{\sigma_1}}   +  \frac{\frac{\sigma_1}{\sigma_2}\phi'\left(\frac{\sigma_1}{\sigma_2}\right)}{1+\phi\left(\frac{\sigma_1}{\sigma_2}\right)}} { \frac{1}{1+\frac{\sigma_1}{\sigma_2}}  +  \frac{\frac{\sigma_2}{\sigma_1}\phi'\left(\frac{\sigma_2}{\sigma_1}\right)}{1+\phi\left(\frac{\sigma_2}{\sigma_1}\right)}} \\
&=&\frac{\left(1+\phi\left(\frac{\sigma_1}{\sigma_2}\right)\right) + \left(1+\frac{\sigma_2}{\sigma_1}\right)\frac{\sigma_1}{\sigma_2}\phi'\left(\frac{\sigma_1}{\sigma_2}\right)}{\left(1+\frac{\sigma_2}{\sigma_1}\right)\left(1+\phi\left(\frac{\sigma_1}{\sigma_2}\right)\right)}
\times \frac{\left(1+\frac{\sigma_1}{\sigma_2}\right) \left(1+\phi\left(\frac{\sigma_2}{\sigma_1}\right)\right)}{ \left(1+\phi\left(\frac{\sigma_2}{\sigma_1}\right)\right)+ \left(1+\frac{\sigma_1}{\sigma_2}\right) \frac{\sigma_2}{\sigma_1}\phi'\left(\frac{\sigma_2}{\sigma_1}\right)} \\
&=&\frac{\left(1+\phi\left(\frac{\sigma_1}{\sigma_2}\right)\right) + \left(1+\frac{\sigma_1}{\sigma_2}\right)\phi'\left(\frac{\sigma_1}{\sigma_2}\right)}{\left(1+\frac{\sigma_2}{\sigma_1}\right)\left(1+\phi\left(\frac{\sigma_1}{\sigma_2}\right)\right)}
\times \frac{\left(1+\frac{\sigma_1}{\sigma_2}\right) \left(1+\phi\left(\frac{\sigma_2}{\sigma_1}\right)\right)}{ \left(1+\phi\left(\frac{\sigma_2}{\sigma_1}\right)\right)+ \left(1+\frac{\sigma_2}{\sigma_1}\right) \phi'\left(\frac{\sigma_2}{\sigma_1}\right)} \\
&=&\frac{\left(1+\phi\left(\frac{\sigma_1}{\sigma_2}\right)\right) + \left(1+\frac{\sigma_1}{\sigma_2}\right)\phi'\left(\frac{\sigma_1}{\sigma_2}\right)}
{ \left(1+\phi\left(\frac{\sigma_2}{\sigma_1}\right)\right)+ \left(1+\frac{\sigma_2}{\sigma_1}\right) \phi'\left(\frac{\sigma_2}{\sigma_1}\right)}
\times \frac{\left(1+\frac{\sigma_1}{\sigma_2}\right) \left(1+\phi\left(\frac{\sigma_2}{\sigma_1}\right)\right)}
{\left(1+\frac{\sigma_2}{\sigma_1}\right)\left(1+\phi\left(\frac{\sigma_1}{\sigma_2}\right)\right)}
 \\
&=&\frac{\left(1+\phi\left(\frac{\sigma_1}{\sigma_2}\right)\right)  \left(1+\phi\left(\frac{\sigma_2}{\sigma_1}\right)\right)
+ \left(1+\frac{\sigma_1}{\sigma_2}\right)\phi'\left(\frac{\sigma_1}{\sigma_2}\right) \left(1+\phi\left(\frac{\sigma_2}{\sigma_1}\right)\right)}
{ \left(1+\phi\left(\frac{\sigma_2}{\sigma_1}\right)\right)\left(1+\phi\left(\frac{\sigma_1}{\sigma_2}\right)\right)
+ \left(1+\frac{\sigma_2}{\sigma_1}\right) \phi'\left(\frac{\sigma_2}{\sigma_1}\right)\left(1+\phi\left(\frac{\sigma_1}{\sigma_2}\right)\right)}
\times
\frac{\left(1+\frac{\sigma_1}{\sigma_2}\right)}{\left(1+\frac{\sigma_2}{\sigma_1}\right)}
 \\
\end{eqnarray*}
The ratio on the right is $\frac{\sigma_1}{\sigma_2}$ and we now show that the ratio on the left is $1$. Clearly the first term in both numerator and denominator of the first ratio is $\left(1+\phi\left(\frac{\sigma_1}{\sigma_2}\right)\right)  \left(1+\phi\left(\frac{\sigma_2}{\sigma_1}\right)\right)$, so the ratio is $1$ if the second term in the numerator is the same as the second term in the denominator, {\em i.e.}, if
$$ \left(1+\frac{\sigma_1}{\sigma_2}\right)\phi'\left(\frac{\sigma_1}{\sigma_2}\right) \left(1+\phi\left(\frac{\sigma_2}{\sigma_1}\right)\right) = \left(1+\frac{\sigma_2}{\sigma_1}\right) \phi'\left(\frac{\sigma_2}{\sigma_1}\right)\left(1+\phi\left(\frac{\sigma_1}{\sigma_2}\right)\right)$$
In other words, what we are required to show is that 
$$ \frac{(1+x)\phi'(x)}{1+\phi(x)}=\frac{\left(1+\frac{1}{x}\right)\phi'\left(\frac{1}{x}\right)}{1+\phi\left(\frac{1}{x}\right)}\,. $$

Note that, from Equations~\eqref{phi} and \eqref{q},
$$ 1+\phi(x)=1+\frac{\nu^2}{8}(x-1)+\frac{\nu}{2}q(x) $$
and combining this with Equation~\eqref{phip},
$$\frac{\phi'(x)}{1+\phi(x)}=\frac{\nu}{4q(x)}\,.$$
Thus
$$\frac{\phi'\left(\frac 1x\right)}{1+\phi\left(\frac 1x\right)}=\frac{\nu}{4q\left(\frac 1x\right)}=\frac{\nu x}{4q(x)}\,,$$
since $xq\left(\frac 1x\right)=q(x)$.
Finally, 
$$ \frac{\left(1+\frac{1}{x}\right)\phi'\left(\frac{1}{x}\right)}{1+\phi\left(\frac{1}{x}\right)} = \frac{\left(1+\frac 1x\right)\nu x}{4q(x)} = \frac{(x+1)\nu}{4q(x)} =\frac{(1+x)\phi'(x)}{1+\phi(x)} $$
as desired.

\qed

Of course, Proposition~\ref{prop:2player} implies Consistency in the Alternative model when $t=2$.

Note also that 
$\frac{\sigma_1}{\sigma_2}=\frac{p_{12}+\frac 12 d_{12}}{p_{21}+\frac 12 d_{21}}$ and, by definition, $\frac{s_1}{s_2}=\frac{w_{12}+\frac 12 t_{12}}{w_{21}+\frac 12 t_{21}}$. So Proposition~\ref{prop:2player} implies that $\frac{p_{12}+\frac 12 d_{12}}{p_{21}+\frac 12 d_{21}}=\frac{w_{12}+\frac 12 t_{12}}{w_{21}+\frac 12 t_{21}}$.

\subsection{Balanced tournaments with three players} \label{sec:3player}

Consider the case $t=3$. As discussed previously, we take $n=1$ and then $s_1=s_{12}+s_{13}$, $s_2=(1-s_{12})+s_{23}$, $s_3=(1-s_{13})+(1-s_{23})$, and $s_1+s_2+s_3=3$. Thinking of the reverse problem, the system of maximum likelihood equations can be written
\begin{eqnarray*}
 s_1 = s_{12}+s_{13} &=& \sigma_1 G_1(s_{12},s_{13}) \\
 s_2 = 1-s_{12}+s_{23} &=& \sigma_2 G_2(s_{12},s_{23}) \\
 s_3 = 2-s_{13}-s_{23} &=& \sigma_3 G_3(s_{13},s_{23}) \\
 T &=& \nu H(s_{12},s_{13},s_{23})
\end{eqnarray*}
which is a linear system of $4$ equations in a vector of $4$ unknowns, $x=(s_{12}, s_{13}, s_{23}, T)$, but with a dependency among the first three equations ($s_1+s_2+s_3=3$), so that, written as $Ax=b$, the rank of $A$ is only $3$.

The equations are
\begin{eqnarray*} 
s_{12}+s_{13} &=& \frac{\sigma_1}{\sigma_1+\sigma_2}+\frac{\sigma_1}{\sigma_1+\sigma_3}
-s_{12}f\left(\frac{\sigma_2}{\sigma_1}\right)
-s_{13}f\left(\frac{\sigma_3}{\sigma_1}\right)
+(1-s_{12})f\left(\frac{\sigma_1}{\sigma_2}\right)
+(1-s_{13})f\left(\frac{\sigma_1}{\sigma_3}\right) \\
1-s_{12}+s_{23} &=& \frac{\sigma_2}{\sigma_2+\sigma_1}+\frac{\sigma_2}{\sigma_2+\sigma_3}
-(1-s_{12})f\left(\frac{\sigma_1}{\sigma_2}\right)
-s_{23}f\left(\frac{\sigma_3}{\sigma_2}\right)
+s_{12}f\left(\frac{\sigma_2}{\sigma_1}\right)
+(1-s_{23})f\left(\frac{\sigma_2}{\sigma_3}\right) \\
2-s_{13}-s_{23} &=& \frac{\sigma_3}{\sigma_3+\sigma_1}+\frac{\sigma_3}{\sigma_3+\sigma_2}
-(1-s_{13})f\left(\frac{\sigma_1}{\sigma_3}\right)
-(1-s_{23})f\left(\frac{\sigma_2}{\sigma_3}\right)
+s_{13}f\left(\frac{\sigma_3}{\sigma_1}\right)
+s_{23}f\left(\frac{\sigma_3}{\sigma_2}\right) \\
T &=& \nu\left[s_{12}g\left(\frac{\sigma_2}{\sigma_1}\right)+s_{13}g\left(\frac{\sigma_3}{\sigma_1}\right)+s_{23}g\left(\frac{\sigma_3}{\sigma_2}\right) \right.\\
&& \left. + (1-s_{12})g\left(\frac{\sigma_1}{\sigma_2}\right) + (1-s_{13})g\left(\frac{\sigma_1}{\sigma_3}\right) +(1-s_{23})g\left(\frac{\sigma_2}{\sigma_3}\right)
\right]
\,, 
\end{eqnarray*}
where 
$$f(x)=\frac{x\phi'\left(x\right)}{1+\phi\left(x\right)}\,,\quad g(x)=\frac{\phi_{\nu}(x)}{1+\phi(x)}\,,$$
which can be written:
{\small
\begin{eqnarray*} 
s_{12}\left[1+f\left(\frac{\sigma_2}{\sigma_1}\right)+f\left(\frac{\sigma_1}{\sigma_2}\right) \right]
+s_{13}\left[1+f\left(\frac{\sigma_3}{\sigma_1}\right)+f\left(\frac{\sigma_1}{\sigma_3}\right)\right] 
&=& \frac{\sigma_1}{\sigma_1+\sigma_2}+\frac{\sigma_1}{\sigma_1+\sigma_3}+f\left(\frac{\sigma_1}{\sigma_2}\right)+f\left(\frac{\sigma_1}{\sigma_3}\right) \\
s_{12}\left[-1-f\left(\frac{\sigma_1}{\sigma_2}\right)-f\left(\frac{\sigma_2}{\sigma_1}\right) \right]
+s_{23}\left[1+f\left(\frac{\sigma_3}{\sigma_2}\right)+f\left(\frac{\sigma_2}{\sigma_3}\right)\right] 
&=& -1 +\frac{\sigma_2}{\sigma_2+\sigma_1}+\frac{\sigma_2}{\sigma_2+\sigma_3}-f\left(\frac{\sigma_1}{\sigma_2}\right)+f\left(\frac{\sigma_2}{\sigma_3}\right) \\
s_{13}\left[-1-f\left(\frac{\sigma_1}{\sigma_3}\right)-f\left(\frac{\sigma_3}{\sigma_1}\right) \right]
+s_{23}\left[-1-f\left(\frac{\sigma_2}{\sigma_3}\right)-f\left(\frac{\sigma_3}{\sigma_2}\right)\right] 
&=& -2 +\frac{\sigma_3}{\sigma_3+\sigma_1}+\frac{\sigma_3}{\sigma_3+\sigma_2}-f\left(\frac{\sigma_1}{\sigma_3}\right)-f\left(\frac{\sigma_2}{\sigma_3}\right) \\
s_{12}\nu \left[-g\left(\frac{\sigma_2}{\sigma_1}\right)+g\left(\frac{\sigma_1}{\sigma_2}\right)\right]+s_{13}\nu\left[-g\left(\frac{\sigma_3}{\sigma_1}\right)+g\left(\frac{\sigma_1}{\sigma_3}\right)\right] &+& s_{23}\nu\left[-g\left(\frac{\sigma_3}{\sigma_2}\right)+g\left(\frac{\sigma_2}{\sigma_3}\right)\right] +T \\
&=&  \nu g\left(\frac{\sigma_1}{\sigma_2}\right) + \nu g\left(\frac{\sigma_1}{\sigma_3}\right)+ \nu g\left(\frac{\sigma_2}{\sigma_3}\right)
\end{eqnarray*}
to put it into the form $Ax=b$.

\begin{prop} The Alternative model is not Consistent when $t=3$, {\em i.e.,} there are balanced designs for which strength parameters do not respect the order of scores.
\end{prop}
{\bf Proof:}
To construct a family of counterexamples, let us first suppose that $\frac{\sigma_1}{\sigma_2}=1$ and $\frac{\sigma_1}{\sigma_3}=5$, which imply that $\frac{\sigma_2}{\sigma_3}=5$, and suppose also that $\nu=4$. We calculate the values needed to construct the coefficients in the linear system as follows.
\begin{eqnarray*}
\frac{\sigma_1}{\sigma_1+\sigma_2}&=&\frac{1}{1+\frac{\sigma_2}{\sigma_1}}=\frac 12\,; \quad \frac{\sigma_1}{\sigma_1+\sigma_3}=\frac{1}{1+\frac{\sigma_3}{\sigma_1}}=\frac 56 \\
\frac{\sigma_2}{\sigma_2+\sigma_1}&=&\frac{1}{1+\frac{\sigma_1}{\sigma_2}}=\frac 12\,; \quad \frac{\sigma_2}{\sigma_2+\sigma_3}=\frac{1}{1+\frac{\sigma_3}{\sigma_2}}=\frac 56 \\
\frac{\sigma_3}{\sigma_3+\sigma_1}&=&\frac{1}{1+\frac{\sigma_1}{\sigma_3}}=\frac 16\,; \quad \frac{\sigma_3}{\sigma_3+\sigma_2}=\frac{1}{1+\frac{\sigma_2}{\sigma_3}}=\frac 16 \\
\end{eqnarray*}
Also, $\phi(x)=2(x-1)+2\sqrt{(x-1)^2+x}$, $\phi'(x)=2+\frac{2x-1}{\sqrt{(x-1)^2+x}}$, and $\phi_{\nu}(x)=(x-1) + \frac{(x-1)^2+\frac{x}{2}  }{\sqrt{(x-1)^2+x}}$, so
\begin{eqnarray*}
\phi\left(\frac 15 \right) &=& -\frac 85 +\frac{2\sqrt{21}}{5} \,,\quad \quad \phi'\left(\frac 15 \right) = 2-\frac{\sqrt{21}}{7} \,,\qquad \phi_{\nu}\left(\frac{1}{5}\right)= -\frac 45 + \frac{37\sqrt{21}}{210} \\
\phi(1) &=&  2 \,,\qquad\qquad\quad \phi'(1)=3\,,\qquad \phi_{\nu}(1)= \frac{1}{2} \\
\phi(5) &=&  8+2\sqrt{21} \,,\qquad \phi'(5)=2+\frac{3\sqrt{21}}{7}\,,\qquad \phi_{\nu}(5)= 4+\frac{37\sqrt{21}}{42}
\end{eqnarray*}
\begin{eqnarray*}
f\left(\frac{\sigma_1}{\sigma_2}\right) = \frac{\frac{\sigma_1}{\sigma_2}\phi'\left(\frac{\sigma_1}{\sigma_2}\right) }{1+\phi\left(\frac{\sigma_1}{\sigma_2}\right)} &=&  1 \,,\qquad f\left(\frac{\sigma_2}{\sigma_1}\right) =\frac{\frac{\sigma_2}{\sigma_1}\phi'\left(\frac{\sigma_2}{\sigma_1}\right)}{1+\phi\left(\frac{\sigma_2}{\sigma_1}\right)}=1 \\
f\left(\frac{\sigma_1}{\sigma_3}\right) = \frac{\frac{\sigma_1}{\sigma_3}\phi'\left(\frac{\sigma_1}{\sigma_3}\right) }{1+\phi\left(\frac{\sigma_1}{\sigma_3}\right)} &=& \frac{5\sqrt{21}}{21} \,,\qquad f\left(\frac{\sigma_3}{\sigma_1}\right) = \frac{\frac{\sigma_3}{\sigma_1}\phi'\left(\frac{\sigma_3}{\sigma_1}\right)}{1+\phi\left(\frac{\sigma_3}{\sigma_1}\right)}=\frac{\sqrt{21}}{21} \\
f\left(\frac{\sigma_2}{\sigma_3}\right) = \frac{\frac{\sigma_2}{\sigma_3}\phi'\left(\frac{\sigma_2}{\sigma_3}\right) }{1+\phi\left(\frac{\sigma_2}{\sigma_3}\right)} &=&  \frac{5\sqrt{21}}{21} \,,\qquad f\left(\frac{\sigma_3}{\sigma_2}\right) =  \frac{\frac{\sigma_3}{\sigma_2}\phi'\left(\frac{\sigma_3}{\sigma_2}\right)}{1+\phi\left(\frac{\sigma_3}{\sigma_2}\right)}=\frac{\sqrt{21}}{21}  \\
\end{eqnarray*}
\begin{eqnarray*}
g\left(\frac{\sigma_1}{\sigma_2}\right) = \frac{\phi_{\nu}\left(\frac{\sigma_1}{\sigma_2}\right) }{1+\phi\left(\frac{\sigma_1}{\sigma_2}\right)} &=&  \frac 16 \,,\qquad g\left(\frac{\sigma_2}{\sigma_1}\right) =\frac{\phi_{\nu}\left(\frac{\sigma_2}{\sigma_1}\right)}{1+\phi\left(\frac{\sigma_2}{\sigma_1}\right)}=\frac 16 \\
g\left(\frac{\sigma_1}{\sigma_3}\right) = \frac{\phi_{\nu}\left(\frac{\sigma_1}{\sigma_3}\right) }{1+\phi\left(\frac{\sigma_1}{\sigma_3}\right)} &=& \frac{1}{3}+\frac{\sqrt{21}}{42} \,,\qquad g\left(\frac{\sigma_3}{\sigma_1}\right) = \frac{\phi_{\nu}\left(\frac{\sigma_3}{\sigma_1}\right)}{1+\phi\left(\frac{\sigma_3}{\sigma_1}\right)}= \frac 13  -\frac{\sqrt{21}}{14}\\
g\left(\frac{\sigma_2}{\sigma_3}\right) = \frac{\phi_{\nu}\left(\frac{\sigma_2}{\sigma_3}\right) }{1+\phi\left(\frac{\sigma_2}{\sigma_3}\right)} &=&  \frac{1}{3}+\frac{\sqrt{21}}{42} \,,\qquad 
g\left(\frac{\sigma_3}{\sigma_2}\right) =  \frac{\phi_{\nu}\left(\frac{\sigma_3}{\sigma_2}\right)}{1+\phi\left(\frac{\sigma_3}{\sigma_2}\right)}=\frac 13  -\frac{\sqrt{21}}{14} \\
\end{eqnarray*}
Finally, the linear system is 
\begin{eqnarray*}
    \left(1+1+1\right)s_{12}+\left(1+\frac{\sqrt{21}}{21}+\frac{5\sqrt{21}}{21} \right)s_{13} &=& 
    \frac 12 +\frac 56 + 1 +  \frac{5\sqrt{21}}{21} \\
    \left(-1-1-1\right)s_{12}+\left(1+\frac{\sqrt{21}}{21}+\frac{5\sqrt{21}}{21} \right)s_{23} &=& -1 +
    \frac 12 +\frac 56 - 1 +  \frac{5\sqrt{21}}{21} \\
    \left(-1-\frac{5\sqrt{21}}{21} -\frac{\sqrt{21}}{21} \right)s_{13}+\left(-1-\frac{5\sqrt{21}}{21} -\frac{\sqrt{21}}{21} \right)s_{23} &=& -2 +
    \frac 16 +\frac 16 - \frac{5\sqrt{21}}{21} -\frac{5\sqrt{21}}{21}  \\
    4\left(\frac 16 - \frac 16 \right)s_{12} + 4\left(\frac 13+\frac{\sqrt{21}}{42} -\frac 13 + \frac{\sqrt{21}}{14}\right)s_{13} + 4\left(\frac 13+\frac{\sqrt{21}}{42} -\frac 13 + \frac{\sqrt{21}}{14}\right)s_{23} &+& \\
    T &=& 4\left(\frac 16 + \frac 13 +\frac{\sqrt{21}}{42} + \frac 13 +\frac{\sqrt{21}}{42}\right)
\end{eqnarray*}
which reduces to $Ax=b$ with
$$ A=\left[\begin{array}{cccc}
3 & \left(1+\frac{2\sqrt{21}}{7}\right)  & 0 & 0 \\
-3 & 0 & \left(1+\frac{2\sqrt{21}}{7}\right) & 0 \\
0 & \left(-1-\frac{2\sqrt{21}}{7}\right) & \left(-1-\frac{2\sqrt{21}}{7}\right) & 0 \\
0 & \left(\frac{8\sqrt{21}}{21}\right) &  \left( \frac{8\sqrt{21}}{21}\right) & 1
\end{array}\right]\,,\quad x=\left[\begin{array}{c} s_{12} \\ s_{13} \\ s_{23} \\T \end{array}\right]\,, \quad b=\left[\begin{array}{c} \left(\frac{7}{3}+\frac{5\sqrt{21}}{21}\right) \\ \left(-\frac 23 +\frac{5\sqrt{21}}{21}\right) \\ \left(-\frac{5}{3} -\frac{10\sqrt{21}}{21}\right) \\ \left(\frac{10}{3} + \frac{4\sqrt{21}}{21}\right) \end{array}\right]  $$
As anticipated, the third equation is dependent on the first two, so we can write the system as
$$ \left[\begin{array}{cccc}
3 & \left(1+\frac{2\sqrt{21}}{7}\right)  & 0 & 0 \\
-3 & 0 & \left(1+\frac{2\sqrt{21}}{7}\right) & 0 \\
0 & \left(\frac{8\sqrt{21}}{21}\right) &  \left( \frac{8\sqrt{21}}{21}\right) & 1
\end{array}\right] \left[\begin{array}{c} s_{12} \\ s_{13} \\ s_{23} \\T \end{array}\right] = \left[\begin{array}{c} \left(\frac{7}{3}+\frac{5\sqrt{21}}{21}\right) \\ \left(-\frac 23 +\frac{5\sqrt{21}}{21}\right) \\ \left(\frac{10}{3} + \frac{4\sqrt{21}}{21}\right) \end{array}\right]  $$
and its solution is
$$ x = \left[\begin{array}{c} s_{12} \\ s_{13} \\ s_{23} \\T \end{array}\right] = \left[\begin{array}{c} \left(\frac{2}{9}-\frac{5\sqrt{21}}{63}\right)\\ \left(\frac 53\right) \\ 0 \\ \left(\frac{10}{3} -\frac{28\sqrt{21}}{63}\right) \end{array}\right] +u\left[\begin{array}{c}\left( \frac 13 +\frac{2\sqrt{21}}{21}\right) \\ -1 \\ 1 \\ 0 \end{array}\right] $$
for arbitrary parameter $u\in [\frac 23,1]$ (to keep each $s_{ij}\in [0,1]$). When $u=\frac 23$, we have solution 
$$ \left[\begin{array}{c} s_{12} \\ s_{13} \\ s_{23} \\T \end{array}\right] = \left[\begin{array}{c} \left(\frac{4}{9}-\frac{\sqrt{21}}{63}\right)\\ 1  \\ \frac 23 \\ \left(\frac{10}{3} -\frac{28\sqrt{21}}{63}\right) \end{array}\right] $$
and when $u=1$, we have solution
$$ \left[\begin{array}{c} s_{12} \\ s_{13} \\ s_{23} \\T \end{array}\right] = \left[\begin{array}{c} \left(\frac{5}{9}+\frac{\sqrt{21}}{63}\right)\\ \frac 23  \\ 1 \\ \left(\frac{10}{3} -\frac{28\sqrt{21}}{63}\right) \end{array}\right]\,. $$
Note that 
$$ s_1=s_{12}+s_{13}=\left(\frac{17}{9} - \frac{5\sqrt{21}}{63}\right)+ \left(-\frac 23 +\frac{2\sqrt{21}}{21}\right)u$$ 
and 
$$ s_2=1-s_{12}+s_{23}=\left(\frac{7}{9}+\frac{5\sqrt{21}}{63}\right)+\left(\frac 23 -\frac{2\sqrt{21}}{21}\right)u\,. $$
Thus,
$$ s_1-s_2 = \left(\frac {10}{9} -\frac{10\sqrt{21}}{63}\right)+\left(-\frac 43 +\frac{4\sqrt{21}}{21}\right)u\,,
$$
which for $u\in[\frac 23,1]$ means
$$ s_1-s_2 \in \left[-\frac 29 +\frac{2\sqrt{21}}{63},\quad \frac 29 - \frac{2\sqrt{21}}{63}\right]\approx [-0.0767,0.0767]\,.$$
Alternatively, 
$$ \frac{s_1}{s_2} \in \left[\frac{251+6\sqrt{21}}{295}, \frac{251-6\sqrt{21}}{211}\right]\approx [0.9441,1.0593]\,.$$
This means that even though $\sigma_1=\sigma_2$, there can be a positive or negative difference between $s_1$ and $s_2$.
If we perturb $\frac{\sigma_1}{\sigma_2}$ slightly from $1$, we can get a choice of scores in the opposite order to the strength parameters.

For example, if we let $\frac{\sigma_1}{\sigma_2}=\frac{24}{25}=0.96$, we can follow the same computations as above to find 
$$ x = \left[\begin{array}{c} s_{12} \\ s_{13} \\ s_{23} \\T \end{array}\right] \approx \left[\begin{array}{c} -23.7321 \\ 32.2866 \\ -30.7852 \\ 0 \end{array}\right] +u\left[\begin{array}{c} 18.8289 \\ -24.4503 \\ 24.5831 \\ 1 \end{array}\right] $$
where we must have $u\in [1.2796, 1.2930]$. 
At $u=1.2796$, we have solution 
\begin{equation} \left[\begin{array}{c} s_{12} \\ s_{13} \\ s_{23} \\T \end{array}\right] = \left[\begin{array}{c} 0.3614  \\ 1  \\ 0.6714 \\ 1.2796 \end{array}\right] \label{eq:solution}
\end{equation}
and at $u=1.2930$, we have solution 
$$ \left[\begin{array}{c} s_{12} \\ s_{13} \\ s_{23} \\T \end{array}\right] = \left[\begin{array}{c} 0.6132  \\ 0.6731  \\ 1 \\ 1.2930 \end{array}\right]\,. $$
Now, 
$$ \frac{s_1}{s_2}=\frac{s_{12}+s_{13}}{1-s_{12}+s_{23}} \in [0.9275,1.03935]\,,$$
and, in particular, if we take $u=1.2796$, we have $\frac{s_1}{s_2}=1.03935$, whereas we took $\frac{\sigma_1}{\sigma_2}=0.96$. So it is possible for scores ordered such that $s_1>s_2$ to nevertheless produce strength estimates such that $\sigma_1<\sigma_2$.

To concretize this last example, we can approximate the scores by actual numbers of wins and draws between players, and calculate the maximum likelihood estimates for the parameters in the usual forwards way. We use $n=100$ now, rather then $n=1$ to allow for integer numbers of wins and draws.

Let 
\begin{equation} (w_{ij}) = \left[\begin{array}{ccc}
0 & 4 & 100 \\
32 & 0 & 35 \\
0 & 1 & 0
\end{array}\right]\,,\quad (t_{ij})=\left[\begin{array}{ccc}
0 & 64 & 0 \\
64 & 0 & 64 \\
0 & 64 & 0
\end{array}\right]\,,\quad (s_{ij})=(w_{ij})+\frac 12 (t_{ij})=\left[\begin{array}{ccc}
0 & 36 & 100 \\
64 & 0 & 67 \\
0 & 33 & 0
\end{array}\right]\,,
\label{counterexample}
\end{equation}
corresponding to~\eqref{eq:solution},
and the total scores for each player are
$$ (s_1,s_2,s_3)=(136, 131, 33)\,.$$
Then, running the iteration~\eqref{iteration_alternative_practical} gives (normalized so that $\sum_i\sigma_i=1$):
$$ \sigma_1 = 0.4456, \quad \sigma_2 = 0.4648, \quad \sigma_3 = 0.0896, \quad \nu=3.9769$$
so that 
$$ \frac{\sigma_1}{\sigma_2}= 0.9586, \quad \frac{\sigma_1}{\sigma_3}=4.9725, \quad \frac{\sigma_2}{\sigma_3}=5.1871\,.$$
Thus, $s_1>s_2$, but $\sigma_1 < \sigma_2$. \qed

\subsection{Balanced tournaments with four or more players}
When $t=4$, there are $t=4$ in independent parameters, and $\frac 12 t(t-1)+1=7$ independent data values in the system of maximum-likelihood equations. Thus, there is now a $3$-dimensional family of data values, all of which predict the same set of parameters, $(\bm{\sigma},\nu)$. This gives even more room to find sets of data values such that scores are not ordered with strength parameters. Following a procedure similar to that of Section~\ref{sec:3player}, we find the following example, again with $n=100$.

Let 
$$ (w_{ij}) = \left[\begin{array}{cccc}
0 & 16 & 6 & 35 \\
1 & 0 & 36 & 53 \\
1 & 1 & 0 & 2 \\
1 & 1 & 1 & 0
\end{array}\right]\,,\quad (t_{ij})=\left[\begin{array}{cccc}
0 & 83 & 93 & 64 \\
83 & 0 & 63 & 46 \\
93 & 63 & 0 & 97 \\
64 & 46 & 97 & 0
\end{array}\right]\,,\quad (s_{ij})=(w_{ij})+\frac 12 (t_{ij})=\left[\begin{array}{cccc}
0 & 57.5 & 52.5 & 67 \\
42.5 & 0 & 67.5 & 76 \\
47.5 & 32.5 & 0 & 50.5 \\
33 & 24 & 49.5 & 0
\end{array}\right]$$
and the total scores for each player are
$$ (s_1,s_2,s_3,s_4)=(177, 186, 130.5,106.5)\,.$$
Then, running the iteration~\eqref{iteration_alternative_practical} gives (normalized so that $\sum_i\sigma_i=1$):
$$ \sigma_1 = 0.3370, \quad \sigma_2 = 0.3188, \quad \sigma_3 = .1873, \quad \sigma_4 = .1570, \quad \nu=16.1409\,.$$
Thus, $s_1<s_2$, but $\sigma_1 > \sigma_2$.

For tournaments with more players, there are even more extra degrees of freedom in the scores, as a result of the fact that the maximum likelihood equations depend on the individual $s_{ij}$, not just the total player scores, $s_i$. One might consider trying a model with more parameters, such as a draw-propensity parameter, $\nu_i$ for each player, and letting $d_{ij}=\sqrt{\nu_i \nu_j}\sqrt{p_{ij}p_{ji}}$, but this only doubles the number of parameters, whereas the number of data values grows quadratically. In order to have enough parameters to match this growth, one might consider a model with a different draw propensity for each pair of players $\nu_{ij}=\nu_{ji}$. However, then the maximum likelihood equations involve the individual $t_{ij}$ as well as $s_{ij}$, so there are still not enough parameters.

\section{Reassessment}\label{sec:reassessment}
It is important to note that when $\nu=2$, the Alternative model is equivalent to the Davidson model by a scaling, as shown in Equations~\eqref{SigmaRatioSimple} and \eqref{pijEdwardsSimple}. Thus, in this case, the Alternative model has no problem with inconsistency of ordering of strength parameters by scores in balanced tournaments. Property~\ref{consistency_property} holds for the Alternative model when $\nu=2$.

In fact, we argue that $\nu=2$ is the correct value to use in estimating player strengths. Davidson's model promises to estimate both strength parameters and draw propensity at the same time from the same set of data. In the Alternative model, we corrected for Davidson's first flawed assumption, by assuming strength ratios depended on score ratios, Property~\ref{alt_strength_ratio}, not just win ratios, Property~\ref{Luce_Condition}, but we still tried to estimate the draw-propensity parameter at the same time as the strength parameters. The unrealistic implication is that strength parameters do not depend only on scores, but also on draw propensity (and thus, numbers of draws that make up the scores). It can reasonably be argued that two draws are equivalent to a win and a loss between a pair of players in assessing relative strength, so that a score of $+2-1=2$ should give the same assessment as $+3-2=0$, even though the draw propensity is clearly different. When $\nu=2$, a pair of draws is treated exactly the same as a win and a loss (one could even say that a draw was treated as half a win). This is essentially what Glickman assumed in his 1999 paper~\cite{Glickman1999}, but now we have a justification for this in the context of a model with a draw-propensity parameter, via Equation~\eqref{pijEdwardsSimple}. Glickman argued that a draw should count as a half-win followed by a half-loss, and so that if $p$ is the probability of a win by the first player, a single draw should contribute a term $\sqrt{p(1-p)}$ to the likelihood function. In fact, it should contribute $2\sqrt{p(1-p}$ to indclude the possibilities of a win followed by a loss (10) and a loss followed by a win (01). In our model, when $\nu=2$, we have $d_{ij}=\frac{2\sigma_i\sigma_j}{(\sigma_i+\sigma_j)^2}=2\sqrt{p_{ij}p_{ji}}$, which is twice the probability of a win and a loss in the Bradley-Terry model.

\subsection{Constrained Alternative (and Constrained Davidson) model}
Fixing $\nu=2$ for the estimation of strength parameters does not in itself allow for estimation of draw propensity, which we still require. But draw propensity should be estimated using meaningful strength parameters.

We argue, therefore, that a more correct procedure is to constrain the Alternative model to $\nu=2$ for the estimation of strength parameters, and thus iterate via
\begin{eqnarray*}
    \sigma_k &\leftarrow& \frac{s_k}{G_k(\bm{\sigma},2)}\,,\quad k=1,\ldots,t, \\
    \nu &\leftarrow& \frac{T}{H(\bm{\sigma}, \nu)}\,.
\end{eqnarray*} 
This procedure effectively estimates strength parameters with $\nu$ fixed at $2$, as it should, and then keeping the resulting strength parameters now fixed, estimates the value of $\nu$ that maximizes the likelihood for this (correct) set of strength parameters.
Note that the two stages of the estimation here use different likelihood functions. The first uses $L(\bm{\sigma},2)$, the likelihood of the observed results when $\nu=2$, to get an estimate $\bm{\sigma}=\bm{\hat{\sigma}}$, and the second stage uses $L(\bm{\hat{\sigma}},\nu)$ to estimate $\nu$. The first stage is exactly the Bradley-Terry model considering each comparison to be two comparisons, such that draws are $10$ (win+loss), while wins are $11$, and losses are $00$. The unconstrained likelihood function $L(\bm{\sigma},\nu)$ allows for too much freedom and sets of parameters that are not sensible. Finally, it is clear that the issue of convergence of the iteration, which we avoided in Section~\ref{sec:alternative_likelihood}, is not a problem when $\nu=2$, since the Bradley-Terry and Davidson methods are known to be convergent.

Of course, the same procedure could be applied to the Davidson model: use $\nu=2$ to estimate strength parameters, $\bm{\pi}$, and then use the resulting strength parameters to estimate $\nu$. Thus, we have four models to evaluate: the unconstrained Davidson model, the Constrained Davidson model (with $\nu=2$ fixed for the strength parameter calculation), the unconstrained Alternative model, and the Constrained Alternative model.

\section{Test cases for the four models} \label{sec:tests}

\subsection{Test 1}
The following test case is instructive.
$$(w_{ij})=\left[\begin{array}{cc} 0 & 1+\epsilon \\ \epsilon & 0\end{array}\right] \,,\quad (t_{ij})=\left[\begin{array}{cc} 0 & 4-2\epsilon \\ 4-2\epsilon & 0\end{array}\right]\quad \implies \quad (s_{ij})=\left[\begin{array}{cc} 0 & 3 \\ 2 & 0\end{array}\right] $$ 
with $\epsilon$ going from $2$ down to $0$.
 
For the {\bf Davidson model}, the maximum log-likelihood equations are
\begin{eqnarray*}
\frac{\partial}{\partial \pi_1}\log(L) &=& \frac{3}{\pi_1}-\frac{5\left(1+\frac{\nu}{2}\sqrt{\frac{\pi_2}{\pi_1}}\right)}{\pi_1+\pi_2+\nu\sqrt{\pi_1\pi_2}} = 0 \\
\frac{\partial}{\partial \pi_2}\log(L) &=& \frac{2}{\pi_2}-\frac{5\left(1+\frac{\nu}{2}\sqrt{\frac{\pi_1}{\pi_2}}\right)}{\pi_1+\pi_2+\nu\sqrt{\pi_1\pi_2}} = 0 \\
\frac{\partial}{\partial \nu}\log(L) &=& \frac{4-2\epsilon}{\nu}-\frac{5\left(\sqrt{\pi_1\pi_2}\right)}{\pi_1+\pi_2+\nu\sqrt{\pi_1\pi_2}} = 0 \\
\end{eqnarray*}
but the first two are necessarily equivalent.
The second and third give:
$$ \frac{\pi_1}{\pi_2}+1+\nu\sqrt{\frac{\pi_1}{\pi_2}}=\frac{5}{2}\left(1+\frac{\nu}{2}\sqrt{\frac{\pi_1}{\pi_2}}\right) = \frac{5}{4-2\epsilon}\nu\sqrt{\frac{\pi_1}{\pi_2}}$$
So
$$\frac 12 =\left(\frac{1}{4-2\epsilon}-\frac{1}{4}\right)\nu\sqrt{\frac{\pi_1}{\pi_2}}$$
or
$$\frac{4-2\epsilon}{\epsilon}=\nu\sqrt{\frac{\pi_1}{\pi_2}}$$
and
$$\frac{\pi_1}{\pi_2}+1+\frac{4-2\epsilon}{\epsilon}=\frac{5}{\epsilon}$$
or
$$ \frac{\pi_1}{\pi_2}=\frac{1+\epsilon}{\epsilon}$$
which $\to \infty$ as $\epsilon\to 0$. Also, then
$$\nu=\frac{4-2\epsilon}{\epsilon}\sqrt{\frac{\epsilon}{1+\epsilon}}=\frac{4-2\epsilon}{\sqrt{\epsilon}\sqrt{1+\epsilon}}$$
which also $\to\infty$ as $\epsilon\to 0$. Finally,
\begin{eqnarray*} 
p_{12} &=& \frac{\frac{\pi_1}{\pi_2}}{\frac{\pi_1}{\pi_2}+1+\nu\sqrt{\frac{\pi_1}{\pi_2}}} = \frac{\frac{1+\epsilon}{\epsilon}}{\frac{1+\epsilon}{\epsilon}+1+\frac{4-2\epsilon}{\epsilon}} = \frac{1+\epsilon}{1+\epsilon+\epsilon+4-2\epsilon}=\frac{1+\epsilon}{5} \\
p_{21} &=& \frac{1}{\frac{\pi_1}{\pi_2}+1+\nu\sqrt{\frac{\pi_1}{\pi_2}}} = \frac{1}{\frac{1+\epsilon}{\epsilon}+1+\frac{4-2\epsilon}{\epsilon}} = \frac{\epsilon}{1+\epsilon+\epsilon+4-2\epsilon}=\frac{\epsilon}{5} \\
d_{12} &=& \frac{\nu\sqrt{\frac{\pi_1}{\pi_2}}}{\frac{\pi_1}{\pi_2}+1+\nu\sqrt{\frac{\pi_1}{\pi_2}}} = \frac{\frac{4-2\epsilon}{\epsilon}}{\frac{1+\epsilon}{\epsilon}+1+\frac{4-2\epsilon}{\epsilon}} = \frac{4-2\epsilon}{1+\epsilon+\epsilon+4-2\epsilon}=\frac{4-2\epsilon}{5}\,. \\
\end{eqnarray*}
Thus, as $\epsilon\to 0$, the probabilities approach $p_{12}=\frac 15$, $p_{21}=0$, and $d_{12}=\frac 45$, corresponding to the scores, $w_{12}=1$, $w_{21}=0$, and $t_{12}=4$ as expected, but at the cost of an infinite separation in player strengths, which does not correspond to the fact that the players' scores were $s_1=3$ and $s_2=2$ (regardless of the value of $\epsilon$ here).

For the {\bf Constrained Davidson model} ($\nu=2$ to estimate strength parameters), the maximum log-likelihood equations are
\begin{eqnarray*}
\frac{\partial}{\partial \pi_1}\log(L) &=& \frac{3}{\pi_1}-\frac{5\left(1+\sqrt{\frac{\pi_2}{\pi_1}}\right)}{\pi_1+\pi_2+2\sqrt{\pi_1\pi_2}} = 0 \\
\frac{\partial}{\partial \pi_2}\log(L) &=& \frac{2}{\pi_2}-\frac{5\left(1+\sqrt{\frac{\pi_1}{\pi_2}}\right)}{\pi_1+\pi_2+2\sqrt{\pi_1\pi_2}} = 0 \\
\frac{\partial}{\partial \nu}\log(L) &=& \frac{4-2\epsilon}{\nu}-\frac{5\left(\sqrt{\pi_1\pi_2}\right)}{\pi_1+\pi_2+\nu\sqrt{\pi_1\pi_2}} = 0 \\
\end{eqnarray*}
but, again, the first two are necessarily equivalent. The first gives
$$ \frac{\pi_1}{\pi_2}+1+2\sqrt{\frac{\pi_1}{\pi_2}}=\frac{5}{3}\left(\frac{\pi_1}{\pi_2}+\sqrt{\frac{\pi_1}{\pi_2}}\right) $$
So
$$\frac 23 \frac{\pi_1}{\pi_2}-\frac 13 \sqrt{\frac{\pi_1}{\pi_2}} = 1 = 0$$
and
$$ \sqrt{\frac{\pi_1}{\pi_2}} = \frac 12\left(\frac 12 \pm \sqrt{\frac 14 + 6}\right) =\frac 32$$
(we must pick the plus sign here) so that 
$$ \frac{\pi_1}{\pi_2}=\frac 94 \,. $$
The third equation gives
$$ \frac{\pi_1}{\pi_2}+1+\nu\sqrt{\frac{\pi_1}{\pi_2}}=\frac{5}{4-2\epsilon}\nu\sqrt{\frac{\pi_1}{\pi_2}} $$
which at $\frac{\pi_1}{\pi_2}=\frac 94$ is
$$ \frac 32\left(\frac{1+2\epsilon}{4-2\epsilon}\right)\nu=\frac 94+1 $$
so
$$ \nu = \frac{13}{6}\left(\frac{4-2\epsilon}{1+2\epsilon}\right) $$
which converges, as $\epsilon\to 0$, to $\frac{26}{3}$.  Here, 
\begin{eqnarray*} 
p_{12} &=& \frac{\frac{9}{4}}{\frac{9}{4}+1+\frac{13}{6}\left(\frac{4-2\epsilon}{1+2\epsilon}\right)\frac 32}=\frac{9}{13\left(1+\frac{4-2\epsilon}{1+2\epsilon}\right)}=\frac{9(1+2\epsilon)}{65} \\
p_{21} &=& \frac{1}{\frac{9}{4}+1+\frac{13}{6}\left(\frac{4-2\epsilon}{1+2\epsilon}\right)\frac 32}=\frac{4}{13\left(1+\frac{4-2\epsilon}{1+2\epsilon}\right)}=\frac{4(1+2\epsilon)}{65} \\
d_{12} &=& \frac{\frac{13}{6}\left(\frac{4-2\epsilon}{1+2\epsilon}\right)\frac 32}{\frac{9}{4}+1+\frac{13}{6}\left(\frac{4-2\epsilon}{1+2\epsilon}\right)\frac 32}=\frac{13\left(\frac{4-2\epsilon}{1+2\epsilon}\right)}{13\left(1+\frac{4-2\epsilon}{1+2\epsilon}\right)}=\frac{13(4-2\epsilon)}{65} 
\end{eqnarray*}
which converge as $\epsilon\to 0$ to $\frac{9}{65}$, $\frac{4}{65}$ and $\frac{52}{65}$, respectively.

For the {\bf Alternative model}, we already know by Proposition~\ref{prop:2player}, that $\frac{\sigma_1}{\sigma_2}=\frac{s_1}{s_2}$, which here is $\frac 32$. Note that this is true for any value of $\nu$.
The third maximum log-likelihood equation is now
\begin{eqnarray*}
\frac{\partial}{\partial \nu}\log(L) &=& \frac{4-2\epsilon}{\nu}-\frac{3\phi_{\nu}\left(\frac{\sigma_2}{\sigma_1}\right)}{1+\phi\left(\frac{\sigma_2}{\sigma_1}\right)} - \frac{2\phi_{\nu}\left(\frac{\sigma_1}{\sigma_2}\right)}{1+\phi\left(\frac{\sigma_1}{\sigma_2}\right)} 
=\frac{4-2\epsilon}{\nu}-\frac{3\phi_{\nu}\left(\frac 23 \right)}{1+\phi\left(\frac 23\right)} - \frac{2\phi_{\nu}\left(\frac 32 \right)}{1+\phi\left(\frac 32 \right)}=0 \\
\end{eqnarray*}
or
$$ (4-2\epsilon) -\nu H(\nu) = (4-2\epsilon) -\nu\left(\frac{3\phi_{\nu}\left(\frac 23 \right)}{1+\phi\left(\frac 23\right)} + \frac{2\phi_{\nu}\left(\frac 32 \right)}{1+\phi\left(\frac 32 \right)}\right)=0$$
The left-hand side function here is a strictly positive function of $\nu$ when $\epsilon=0$, so there is no maximum, but has a zero when $0<\epsilon\le 2$. The proof of this requires only the algebra to show that 
\begin{equation}
\nu H(\nu) = \nu\left(\frac{3\phi_{\nu}\left(\frac 23 \right)}{1+\phi\left(\frac 23\right)} + \frac{2\phi_{\nu}\left(\frac 32 \right)}{1+\phi\left(\frac 32 \right)}\right)= \frac{\frac{5}{4}\nu^2+\frac{\nu}{4}\sqrt{\nu^2+96}}{1+\frac{13}{48}\nu^2+\frac{5}{48}\nu\sqrt{\nu^2+96}}\,, \label{nuHnu}
\end{equation}
which is $0$ at $\nu=0$ with $\lim_{\nu\to\infty}\nu H(\nu) = 4 $,
and to show that 
$$ \frac{\partial}{\partial\nu} (\nu H(\nu)) = \frac{\partial}{\partial\nu} \left(\frac{\frac{5}{4}\nu^2+\frac{\nu}{4}\sqrt{\nu^2+96}}{1+\frac{13}{48}\nu^2+\frac{5}{48}\nu\sqrt{\nu^2+96}}\right)=\frac{24+\frac{13}{2}\nu^2+\frac{5}{2}\nu\sqrt{\nu^2+96}}{\left[1+\frac{13}{48}\nu^2+\frac{5}{48}\nu\sqrt{\nu^2+96}\right]^2 \sqrt{\nu^2+96}}>0\,.  $$
Thus, our estimate of $\nu\to\infty$ as $\epsilon\to 0$. In fact, the maximum likelihood estimate is 
\begin{equation}
    \nu = \frac{4-2\epsilon}{\sqrt{\epsilon(\epsilon+1)}}\,, \label{nu_Alternative}
\end{equation}
which can be verified by inserting~\eqref{nu_Alternative} into Equation~\eqref{nuHnu} to get $\nu H(\nu)=4-2\epsilon$.

Now, with this value of $\nu$,
\begin{eqnarray*}
p_{12} &=& \left(\frac{1}{1+\frac{\sigma_2}{\sigma_1}}\right)\left(\frac{1}{1+\phi\left(\frac{\sigma_2}{\sigma_1}\right)}\right)=\left(\frac{1}{1+\frac 23}\right)\left(\frac{1}{1+\phi\left(\frac 23\right)}\right)=\left(\frac 35\right)\left(\frac{1}{1+\frac{\nu}{24}\left(-\nu+\sqrt{\nu^2+96}\right)}\right)=\frac{\epsilon+1}{5} \\
p_{21} &=& \left(\frac{1}{1+\frac{\sigma_1}{\sigma_2}}\right)\left(\frac{1}{1+\phi\left(\frac{\sigma_1}{\sigma_2}\right)}\right)=\left(\frac{1}{1+\frac 32}\right)\left(\frac{1}{1+\phi\left(\frac 32\right)}\right)=\left(\frac 25\right)\left(\frac{1}{1+\frac{\nu}{16}\left(\nu+\sqrt{\nu^2+96}\right)}\right)=\frac{\epsilon}{5} \\
d_{12} &=& 1-p_{12}-p_{21} = \frac{4-2\epsilon}{5}\,. \\
\end{eqnarray*}

For the {\bf Constrained Alternative Model},  we still have $\frac{\sigma_1}{\sigma_2}=\frac{s_1}{s_2}=\frac 32$, since this is true for any $\nu$, including $\nu=2$. Thus, we again estimate $\nu$, $p_{12}$, $p_{21}$, and $d_{12}$ as for the unconstrained Alternative Model.

\begin{table}
\renewcommand{\arraystretch}{1.5}
\begin{centering}
\begin{tabular}{|l|c|c|c|c|c|}
\hline
    Model & $\frac{\pi_1}{\pi_2}$ or $\frac{\sigma_1}{\sigma_2}$ & $\nu$ & $p_{12}$ & $p_{21}$ & $d_{12}$ \\ [0.8ex] \hline\hline
    Davidson & $\frac{1+\epsilon}{\epsilon}$ & $\frac{4-2\epsilon}{\sqrt{\epsilon(1+\epsilon})}$ & $\frac{1+\epsilon}{5}$ & $\frac{\epsilon}{5}$ & $\frac{4-2\epsilon}{5}$ \\ [1.5ex] \hline
    Constrained Davidson & $\frac{9}{4}$ & $\frac{13}{6}\left(\frac{4-2\epsilon}{1+2\epsilon}\right)$ & $\frac{9(1+2\epsilon)}{65}$ & $\frac{4(1+2\epsilon)}{65}$ & $\frac{4-2\epsilon}{5}$ \\ [1.5ex] \hline
    Alternative & $\frac{3}{2}$ & $\frac{4-2\epsilon}{\sqrt{\epsilon(1+\epsilon)}}$ & $\frac{1+\epsilon}{5}$ & $\frac{\epsilon}{5}$ & $\frac{4-2\epsilon}{5}$ \\ [1.5ex] \hline
    Constrained Alternative & $\frac{3}{2}$ & $\frac{4-2\epsilon}{\sqrt{\epsilon(1+\epsilon)}}$ & $\frac{1+\epsilon}{5}$ & $\frac{\epsilon}{5}$ & $\frac{4-2\epsilon}{5}$ \\ [1.5ex] \hline
\end{tabular}
\caption{Estimates of the four models for Test Case 1 where $\frac{s_1}{s_2}=\frac 32$.\label{test1results}}
\end{centering}
\end{table}

We summarize the above results for this example, using each of the four models, in Table~\ref{test1results}.
Note that when $\epsilon=2$ (the case where there were no draws), all four models give $\nu=0$ and $d_{12}=0$, and all but the Constrained Davidson model give $\frac{\pi_1}{\pi_2}=\frac 32$ or $\frac{\sigma_1}{\sigma_2}=\frac 32$, but the Constrained Davidson model finds $\frac{\pi_1}{\pi_2}=\frac 94$, which is the square of the estimate in the Alternative models, consistent with Equations~\eqref{Luce} and \eqref{SigmaRatioSimple}.

As $\epsilon\to 0$, which means that the maximum number of draws ($T=4$) occur, and the second player scores no wins ($s_{21}=0$), the Davidson model tries to make $\frac{\pi_1}{\pi_2}$ as large as possible in order to drive the number of draws to its maximum, and $\lim_{\epsilon\to 0}\frac{\pi_1}{\pi_2}=\infty$. The Constrained Davidson model avoids this, but has the problem that estimation of $\nu$ is independent of the $s_i$ or $s_{ij}$, since, once the strength parameters are fixed by solving the first two maximum likelihood equations, the third depends only on $T$ and $n$. Because only $T=4$ of $n=5$ games are draws, it settles for a finite $\nu=\frac{26}{3}$. Larger $\nu$ make the likelihood of 4 draws out of 5 smaller.  In all the other models, the estimation of $\nu$ depends, directly or indirectly, on the scores. Even in the unconstrained Davidson model, estimation of $\nu$ depends indirectly on $s_{ij}$ through the simultaneous estimation of the $\pi_i$. In the Alternative model, estimation of $\nu$ depends not only on $n$ but directly on the individual $s_{ij}$.
With the extra information used by the other three models that the scores are $s_1=s_{12}=3$ and $s_2=s_{21}=2$, it is clear that the maximum possible number of draws is $T=4$, and thus, the larger the value of $\nu$, the more likely it becomes to achieve this maximum, and the estimate of $\nu$ goes to $\infty$. If one feels that the corresponding $p_{21}=0$ is too low for $0$ wins and $4$ draws out of $5$ games, remember that the same conclusion must apply if all scores are scaled up, and the second player got $0$ wins and $400$ draws out of $500$ games.

Finally, when $\epsilon=\frac 45$, in all models we have $\nu=2$, the case where the actual number of draws corresponds to the number we would expect if each game were two comparisons, and a draw were a $10$ or $01$ result. The Davidson model (constrained and unconstrained are now the same) has $\frac{\pi_1}{\pi_2}=\frac 94$, while the alternative model (constrained and unconstrained) has $\frac{\sigma_1}{\sigma_2}=\frac 32$. All four models estimate $p_{12}=\frac{9}{25}$, $p_{21}=\frac{4}{25}$, and $d_{12}=\frac{12}{25}$.

The fact that the player strength ratio does not reflect the score ratio seems a serious drawback of the Davidson model, extreme in this test case when one player has no wins at all against the other. The Constrained Davidson model has the serious drawback that it settles for a finite $\nu$ when the maximum possible number of draws occurs, preventing the number of draws from going still higher, which in any case is impossible. Correspondingly, it settles on $p_{21}=\frac{4}{65}$ in this case, instead of $0$.

\subsection{Test 2}
Consider the following data:
$$(w_{ij})=\left[\begin{array}{ccc} 0 & 3 & 1+\epsilon \\ 2 & 0 & 0 \\ \epsilon & 0 & 0\end{array}\right] \,,\quad (t_{ij})=\left[\begin{array}{ccc} 0 & 0 & 4-2\epsilon \\ 0 & 0 & 0 \\ 4-2\epsilon & 0 & 0 \end{array}\right]\quad \implies \quad (s_{ij})=\left[\begin{array}{ccc} 0 & 3 & 3 \\ 2 & 0 & 0 \\ 2 & 0 & 0 \end{array}\right] $$ 
in which $w_{12}=3,w_{21}=2,t_{12}=0$ (like Test 1 with $\epsilon=2$), and $w_{13}=1+\epsilon,w_{31}=\epsilon,t_{13}=4-2\epsilon$ (like Test 1 with arbitrary $\epsilon\in [0,2]$), and players two and three do not play against each other. The first player's score against each of the others is $s_{12}=3$ and $s_{13}=3$, while the other player's scores are $s_{21}=2$ and $s_{31}=2$. Although these pairings were treated very differently by the Davidson model in isolation (Test 1), where $\frac{\pi_1}{\pi_2}=\frac 32$ for $w_{12}=3,w_{21}=2,t_{12}=0$, but $\frac{\pi_1}{\pi_3}=\frac{1+\epsilon}{\epsilon}\to\infty$ as $\epsilon\to 0$ for $w_{13}=1+\epsilon,w_{31}=\epsilon,t_{13}=4-2\epsilon$, here they are all part of the same data set, and the strengths of the second and third players are equal by all four models. However, the strength ratios in the Davidson model still depend on the number of draws. Table~\ref{test2results} summarizes the parameter estimates. The calculations are done in a way that parallels thosse for Test 1, and the algebra is here omitted.

\begin{table}
\renewcommand{\arraystretch}{1.5}
\begin{centering}
\begin{tabular}{|l|c|c|c|c|c|c|}
\hline
    Model & $\frac{\pi_1}{\pi_2}$ or $\frac{\sigma_1}{\sigma_2}$ & $\frac{\pi_1}{\pi_3}$ or $\frac{\sigma_1}{\sigma_3}$ & $\nu$ & $p_{12}=p_{13}$ & $p_{21}=p_{31}$ & $d_{12}=d_{13}$  \\ [0.8ex] \hline\hline
    Davidson & $\frac{4+\epsilon}{2+\epsilon}$ & $\frac{4+\epsilon}{2+\epsilon}$ & $\frac{4-2\epsilon}{\sqrt{(2+\epsilon)(4+\epsilon)}}$ & $\frac{4+\epsilon}{10}$ & $\frac{2+\epsilon}{10}$ & $\frac{2-\epsilon}{5}$ \\ [1.5ex] \hline
    Constrained Davidson & $\frac{9}{4}$ & $\frac{9}{4}$ & $\frac{13}{6}\frac{(2-\epsilon)}{(3+\epsilon)}$ & $\frac{9(3+\epsilon)}{65}$ & $\frac{4(3+\epsilon)}{65}$ & $\frac{2-\epsilon}{5}$  \\ [1.5ex] \hline
    Alternative & $\frac{3}{2}$ & $\frac{3}{2}$ & $\frac{4-2\epsilon}{\sqrt{(2+\epsilon)(4+\epsilon)}}$ & $\frac{4+\epsilon}{10}$ & $\frac{2+\epsilon}{10}$ & $\frac{2-\epsilon}{5}$  \\ [1.5ex] \hline
    Constrained Alternative & $\frac{3}{2}$ & $\frac{3}{2}$ & $\frac{4-2\epsilon}{\sqrt{(2+\epsilon)(4+\epsilon)}}$ & $\frac{4+\epsilon}{10}$ & $\frac{2+\epsilon}{10}$ & $\frac{2-\epsilon}{5}$  \\ [1.5ex]  \hline
\end{tabular}
\caption{Estimates of the four models for Test Case 2, where $\frac{s_1}{s_2}=\frac{s_1}{s_3}=\frac 32$.\label{test2results}}
\end{centering}
\end{table}

\subsection{Test 3}
Consider the data in Equations~\eqref{counterexample}. Table~\ref{test3results} summarizes the main estimates by the four models. Note the high player strength ratios of the Davidson model, even in relation to $\left(\frac{s_1}{s_2}\right)^2$ and the fact already observed, that the unconstrained Alternative model estimates $\frac{\sigma_1}{\sigma_2}<1$ even though $\frac{s_1}{s_2}>1$. 
\begin{table}
\renewcommand{\arraystretch}{1.5}
\begin{centering}
\begin{tabular}{|l|c|c|c|}
\hline
    Model & $\frac{\pi_1}{\pi_2}$ or $\frac{\sigma_1}{\sigma_2}$ & $\frac{\pi_1}{\pi_3}$ or $\frac{\sigma_1}{\sigma_3}$ & $\nu$  \\ [0.8ex] \hline\hline
    Davidson & $1.24249$ & $78.88099$ & $3.87200$  \\ [1.5ex] \hline
    Constrained Davidson & $1.16984$ & $27.72784$ & $2.81065$  \\ [1.5ex] \hline
    Alternative & $0.95862$ & $4.97246$ & $3.97690$  \\ [1.5ex] \hline
    Constrained Alternative & $1.08159$ & $5.26572$ & $3.63972$  \\ [1.5ex]  \hline
\end{tabular}
\caption{Estimates of the four models for Test Case 3, where $\frac{s_1}{s_2}=\frac{136}{131}$, and $\frac{s_1}{s_3}=\frac{136}{33}$. \label{test3results}}
\end{centering}
\end{table}

\subsection{Conclusions from the tests}
The first three models all have undesirable behaviour. The Davidson model does not handle the first test well, where one player has $0$ wins. In Test 2, when this situation is not isolated, the Davidson model at least finds the second and third players to be of equal strength, but the strength gap with the first player still depends on how much of the scores come from draws. The Constrained Davidson model does not handle Test 1 well either, for a different reason, since when estimating $\nu$, it loses information about individual scores. Both the unconstrained and Constrained Davidson models produce larger strength separations than the Alternative models, where they are closer to score ratios. The unconstrained Alternative model, however, sometimes inverts the ordering of strength parameters in relation to the ordering of scores in balanced tournaments, as in Test 3. 

All of these flaws result from what we argue are faulty assumptions in the first three models. The Davidson model starts with the assumption that player strength ratios should match the ratio of probabilities of wins only, instead of the expected ratio of scores (based on the probability of a win plus half the probability of a draw). Although in the end, the estimation procedure uses only scores and not just wins, the underlying assumption influences results. The second flaw of the Davidson model, shared by the unconstrained Alternative model, is allowing the propensity for draws indicated by the data to influence strength parameter estimation, rather than just using raw scores. The constrained models just use raw scores to estimate strength parameters, guaranteeing that two draws count the same as one win each (or that a draw counts as half a win).

The Constrained Alternative model avoids the flawed assumptions of the other models and produces estimates that make most sense, consistently across the possible data sets, including extreme cases.



We wish to highlight the observation that a pair of draws should be equivalent to a win and a loss in regards to strength estimation. A win and a loss indicates more variance in results than two draws (playing riskier; or more randomness in judging two items, as opposed to more reluctance to assert an unclear difference), but should not change estimates of the strength parameters themselves. So $+1-0=4$ should give the same strength assessment as $+2-1=2$ or $+3-2=0$. And any reasonable model should give equality to players B and C if A against B gets $+1-0=4$ while player A against player C gets $+3-2=0$. All four models do this, at least when all three players are in one tournament, but the unconstrained Davidson model does not when there are two separate matches $A$ versus $B$ and $A$ versus $C$.


Conceptually one way to ensure that a pair of draws is equivalent to a win and a loss is to think of a game of chess as effectively two paired comparison experiments, where if the game outcome is a draw it is modelled as a win and a loss (on the half-game scale), while a won game is modelled as two wins and a lost game as two losses. (One can speculate about two phases of a game, where obtaining an advantage is phase one, and then converting it to a win is phase 2, but any such analogy is liable to fail to cover all types of games.) The linear scaling of playing strengths ($\log(\bm{\sigma})$) is then halved in comparison to Davidson, but otherwise equivalent. Glickman's 1999 model essentially follows this idea, and is shown to be equivalent to assuming $\nu=2$ in estimating playing strengths in our Alternative model. However, our full model still allows estimation of draw propensity, which Glickman's 1999 paper does not. Also, the estimation of playing strengths using $\nu=2$ in either the Davidson or Alternate model is essentially equivalent, up to a scaling, to using the standard Bradley-Terry model with scores allowing half points.

\section{Estimating the number of draws} \label{sec:draws}

We now return to the problem of estimating strength or preference parameters when only the winners of three-way comparisons are known, such as pools in chess. In Section~\ref{sec:pools} we estimated the number of encounters won and played between each pair of players, allowing for the possibility of unplayed encounters and tied pools. Now we have a model to deal with drawn games, we return to the problem of estimating the number of draws that occurred between each pair of players, and using this to estimate strength parameters. In order to do this, we need to determine a value of the draw-propensity parameter, $\nu$, from similar data where information on numbers of draws is known.

From Equation~\eqref{encounterbygame} we have
\begin{eqnarray*} 
e_{BC}=\frac{p_{BC}}{1-d_{BC}}\,. \label{epd}
\end{eqnarray*}
From the available data on the three-way comparisons, we have 
\begin{equation}
    e_{BC}=\frac{s^e_{BC}}{n^e_{BC}}\,. \label{eq:esn}
\end{equation}

Effectively, the Bayesian problem is to estimate parameters from data, but here we have only partial data. This would make it hopeless, except that the only thing missing is the number of draws. Estimating $\nu$ from other data allows us to estimate numbers of draws, and then carry out the normal Bayesian maximum likelihood calculation to estimate strength parameters.

What we have from above in~\eqref{eq:ne} and \eqref{eq:se} are estimates of numbers and scores of encounters, and thus of $e_{BC}$, $e_{BD}$, and $e_{CD}$. What we want now is to estimate $t_{BC}$, $t_{BD}$, and $t_{CD}$, the expected numbers of draws that occurred. To get these, we need estimates of $d_{BC}$, $d_{BD}$, and  $d_{CD}$. Since $d_{BC}=\nu\sqrt{p_{BC}p_{CB}}$, we have
$$e_{BC}=\frac{p_{BC}}{1-\nu\sqrt{p_{BC}p_{CB}}}\,,\quad\mbox{or}\quad e_{BC}(1-\nu\sqrt{p_{BC}p_{CB}})=p_{BC}\,.$$
Also, since $p_{BC}+p_{CB}+d_{BC}=p_{BC}+p_{CB}+\nu\sqrt{p_{BC}p_{CB}}=1$, we can solve for $P_{CB}$ in terms of $P_{BC}$:
$$ p_{CB} + (\nu\sqrt{p_{BC}})\sqrt{p_{CB}} - (1-p_{BC})=0 $$
so
$$ \sqrt{p_{CB}}=\frac{-\nu\sqrt{p_{BC}}+\sqrt{\nu^2 p_{BC} + 4(1-p_{BC})}}{2} = -\frac{\nu}{2}\sqrt{p_{BC}} + \sqrt{\left(\frac{\nu^2}{4}-1\right)p_{BC}+1}$$
and
$$e_{BC}\left(1-\nu\sqrt{p_{BC}}\left[-\frac{\nu}{2}\sqrt{p_{BC}} + \sqrt{\left(\frac{\nu^2}{4}-1\right)p_{BC}+1}\right]\right)=p_{BC}$$
or
$$ e_{BC}\left(1+\frac{\nu^2}{2}p_{BC}-\nu\sqrt{\left(\frac{\nu^2}{4}-1\right)p_{BC}^2+p_{BC}}\right)=p_{BC}\,.$$
Solving this for $p_{BC}$ gives
$$ p_{BC}=\frac{e_{BC}}{1+\nu\sqrt{e_{BC}-e_{BC}^2}}$$
and then from Equation~\eqref{epd},
\begin{equation} 
d_{BC}=1-\frac{p_{BC}}{e_{BC}} = 1-\frac{1}{1+\nu\sqrt{e_{BC}-e_{BC}^2}}=\frac{\nu\sqrt{e_{BC}-e_{BC}^2}}{1+\nu\sqrt{e_{BC}-e_{BC}^2}}\,. \label{eq:dbc}
\end{equation}
Now, given that $B$ eventually wins against $C$, the probability of this happening with $k$ preceding draws is 
$$\frac{d_{BC}^k p_{BC}}{e_{BC}}\,,$$
and the expected number of draws is
$$ \frac{\sum_{k=0}^{\infty}kd_{BC}^k p_{BC}}{e_{BC}} = \frac{p_{BC}}{e_{BC}}\frac{d_{BC}}{(1-d_{BC})^2} = \frac{d_{BC}}{1-d_{BC}}\,.$$
Thus, for every encounter in which $B$ beats $C$, we expect $\frac{d_{BC}}{1-d_{BC}}$ draws to preceed it, and the expected total number of draws in such encounters is $s^e_{BC}\frac{d_{BC}}{1-d_{BC}}$. Similarly, the expected total number of draws in encounters when $C$ beats $B$ is $s^e_{CB}\frac{d_{BC}}{1-d_{BC}}$, and in all, the expected total number of draws between $B$ and $C$ is
\begin{equation}
t_{BC}=s^e_{BC}\frac{d_{BC}}{1-d_{BC}}+s^e_{CB}\frac{d_{BC}}{1-d_{BC}}=n^e_{BC}\frac{d_{BC}}{1-d_{BC}}\,, \label{eq:tnd}
\end{equation}
where $n^e_{BC}$ is the estimated number of encounters between $B$ and $C$, and the estimated number of encounters won by $B$ over $C$ is $s^e_{BC}$, which is also the expected number of won games by $B$ over $C$, and similarly for the other pairs. Thus the expected total number of wins in games by $B$ over $C$ is
$$ w_{BC}=s^e_{BC}\,,\quad w_{CB}=n^e_{BC}-s^e_{BC}\,,$$
and the expected total number of games between $B$ and $C$ is 
$$ n_{BC}=t_{BC}+n^e_{BC}=n^e_{BC}\left(\frac{d_{BC}}{1-d_{BC}}\right)+n^e_{BC}=\frac{n^e_{BC}}{1-d_{BC}}$$
and similarly for the other pairs. 

From Equations~\eqref{eq:dbc} and \eqref{eq:esn}, we have
$$ \frac{d_{BC}}{1-d_{BC}}=\nu\sqrt{e_{BC}-e_{BC}^2} = \nu\sqrt{e_{BC}}\sqrt{1-e_{BC}}=\nu \sqrt{\frac{s^e_{BC}}{n^e_{BC}}}\sqrt{1-\frac{s^e_{BC}}{n^e_{BC}}}=\nu\frac{\sqrt{s^e_{BC}(n^e_{BC}-s^e_{BC})}}{n^e_{BC}}
$$
and, from Equation~\eqref{eq:tnd},
$$ t_{BC} = \nu\sqrt{s^e_{BC}(n^e_{BC}-s^e_{BC})} =\nu\sqrt{s^e_{BC}s^e_{CB}}\,.$$

Finally, to estimate playing strengths from the incomplete data of results of the rounds, using the Constrained Alternative model, we can now use expected numbers of wins and draws between each pair of players:
\begin{equation} \label{eq:wt}
    w_{ij}=s^e_{ij}\,,\quad t_{ij}=\nu\sqrt{s^e_{ij}s^e_{ji}}
\end{equation}
taking $n^e_{ij}$ and $s^e_{ij}$ from~\eqref{eq:ne} and \eqref{eq:se}, and $\nu$ estimated from other data.

\section{Example: Paris 1821} \label{sec:Paris}

In order to get estimates for our historical example from the models with draws, we need to find a  reasonable value of $\nu$. We cannot estimated it from the results of the pools themselves, since the results we have do not contain information on draws at all. Also, the frequency of draws has changed at top levels over the course of history, and is certainly dependent to some degree on the strengths of the players. One possibility is to take the long series of games played in 1834 between de la Bourdonnais and Alexander McDonnell, for which we have precise results, as an indication of the draw rates for players of this calibre in the same era. In this series, de la Bourdonnais won 45 out of 85 games, McDonnell won 27, and 13 were drawn~\cite{Utterberg2005}. Thus, our estimate for the probability of a draw is $\frac{13}{85}$ and the estimated probabilities of wins by the two players are $\frac{45}{85}$ and $\frac{27}{85}$, respectively. Our estimate from~(\ref{DavidsonDraw}) is therefore $\nu=\frac{d_{ij}}{\sqrt{p_{ij}p_{ji}}}=\frac{\frac{13}{85}}{\sqrt{\frac{45}{85}\frac{27}{85}}}=\frac{13}{9\sqrt{15}}\approx 0.372954$. 

However, we can do better. I have compiled a list of match results among seven strong players over the period 1821-1836 for which we have precise information, although two of the players (Deschapelles and Lewis) give different handicaps, and we should consider them as different players when giving different handicaps (they clearly play with different strengths), so we can consider it a set of results over 9 players. The match results are given in Table~\ref{table:9players}.

\begin{table}
\begin{tabular}{llccccll}
 Player & Opponent & Year & Wins & Losses & Draws & Handicap & Refs.\\ \hline
 Deschapelles & Lewis & 1821 & $0$ & $1$ & $2$ & P \& 1 & \cite[p.27]{Sergeant1934}\,\cite[p.878]{Murray1913} \\ 
 de la Bourdonnais & Lewis & 1823 & $5$ & $2$ & $0$ & & \cite[p.882]{Murray1913}\,\cite[p.184]{Hooper1984}\\ 
 Lewis & Walker & 1829 & $1$ & $1$ & $1$ & Knight & \cite[p.190]{Murray1906}\\ 
 McDonnell & Fraser & 1831 &  $3$ & $1$ & $1$ & & \cite[p.35]{Sergeant1934}\\
 de la Bourdonnais & McDonnell & 1834 & $45$ & $27$ & $13$ & & \cite{Utterberg2005}\,\cite[p.47]{Hooper1984} \\ 
 McDonnell & Walker & 1835 &  $10$ & $1$ & $3$ & & \cite[p.35]{Sergeant1934} \\
 Deschapelles & Saint-Amant & 1836 & $1$ & $1$ & $1$ & P \& 2  & \cite[p.88]{Hooper1984}\\ 
 Saint-Amant & Walker & 1836 &  $5$ & $3$ & $1$ & & \cite[p.42]{Sergeant1934} \\
 Saint-Amant & Fraser & 1836 &  $1$ & $0$ & $2$ & & \cite[p.42]{Sergeant1934} \\
\end{tabular}
\caption{Match results of strong players 1821-1836 (P \& 1 = Pawn and Move; P \& 2 = Pawn and Two Moves)\label{table:9players}}
\end{table}
\vspace{5mm}

The matrices of wins and draws are then:
\begin{eqnarray*}
    (w_{ij}) = \left[\begin{array}{ccccccccc}
    0 & 0 & 0 & 0 & 0 & 0 & 0 & 0 & 0 \\
    0 & 0 & 0 & 0 & 0 & 0 & 1 & 0 & 0 \\
    1 & 0 & 0 & 0 & 2 & 0 & 0 & 0 & 0 \\
    0 & 0 & 0 & 0 & 0 & 1 & 0 & 0 & 0 \\
    0 & 0 & 5 & 0 & 0 & 0 & 0 & 45 & 0 \\
    0 & 0 & 0 & 1 & 0 & 0 & 3 & 1 & 0 \\
    0 & 1 & 0 & 0 & 0 & 5 & 0 & 0 & 1 \\
    0 & 0 & 0 & 0 & 27 & 10 & 0 & 0 & 3 \\
    0 & 0 & 0 & 0 & 0 & 0 & 0 & 1 & 0 
    \end{array}\right]\,,\quad (d_{ij}) = \left[\begin{array}{ccccccccc}
    0 & 0 & 2 & 0 & 0 & 0 & 0 & 0 & 0 \\
    0 & 0 & 0 & 0 & 0 & 0 & 1 & 0 & 0 \\
    2 & 0 & 0 & 0 & 0 & 0 & 0 & 0 & 0 \\
    0 & 0 & 0 & 0 & 0 & 1 & 0 & 0 & 0 \\
    0 & 0 & 0 & 0 & 0 & 0 & 0 & 13 & 0 \\
    0 & 0 & 0 & 1 & 0 & 0 & 1 & 3 & 0 \\
    0 & 1 & 0 & 0 & 0 & 1 & 0 & 0 & 2 \\
    0 & 0 & 0 & 0 & 13 & 3 & 0 & 0 & 1 \\
    0 & 0 & 0 & 0 & 0 & 0 & 2 & 1 & 0 
    \end{array}\right]
\end{eqnarray*}
with players in the following order: Deschapelles (giving Pawn and Move odds), Deschapelles (given Pawn and Two Move odds), Lewis, Lewis (giving Knight odds), de la Bourdonnais, Walker, de Saint-Amant, McDonnell, and Fraser. Although these results span 16 years and some of the players may have varied in strength over that period, here we are mainly concerned with the estimation of the draw-propensity parameter, $\nu$.

Running the Constrained Alternative model~\eqref{iteration_alternative_practical} on this gives an estimate of 
\begin{equation}
\nu=0.4814882\,. \label{nufromdata}
\end{equation}
(The unconstrained Davidson model and the unconstrained Alternative model give similar estimates: $\nu=0.4814241$ and $0.4814897$ respectively, but the Constrained Davidson model, as often seems to be the case, gives a somewhat different estimate: $\nu=0.532327$.) 

Taking a reasonable value of the draw propensity among top players in this era to be $\nu=0.4814882$, we now estimate the likely numbers of draws in the series of pools played by de la Bourdonnais, Cochrane, and Deschapelles in 1821.

From~\eqref{eq:wt}, we have that $t_{ij}=\nu\sqrt{s^e_{ij}s^e_{ji}}$, where $n^e_{ij}$ and $s^e_{ij}$ are given in~\eqref{eq:ne} and \eqref{eq:se}, and $\nu$ in~\eqref{nufromdata}. This gives 
\begin{equation}
    t_{BC}=3.57836\,,\quad t_{BD}=5.44349\,,\quad t_{CD}=3.26743\,.
\end{equation}

We also have the estimated numbers of wins, from~\eqref{eq:wt} and \eqref{eq:se}, so
we can use the Constrained Alternative model to estimate playing strength parameters for the three players, using the data
$$
(w_{ij})=\left[\begin{array}{ccc}
0 & 17.54013 & 16.06749 \\
3.14894 & 0 & 4.08641 \\
7.95492 & 11.26933 & 0
\end{array}\right]\,,\quad 
(t_{ij}) = \left[\begin{array}{ccc}
0 & 3.57836 & 5.44349 \\
3.57836 & 0 & 3.26743 \\
5.44349 & 3.26743 & 0
\end{array}\right]\,,
$$
so that
$$
(s_{ij})=(w_{ij})+\frac 12 (t_{ij}) = \left[\begin{array}{ccc}
0 & 19.32931 & 18.78924 \\
4.93812 & 0 & 5.72012 \\
10.67667 & 12.90304 & 0
\end{array}\right]\,.
$$
The resulting strength parameters, normalized so that $\sum_i \sigma_i=1$, are
$$ \bm{\sigma} = ( 0.54821, 0.13929, 0.31250 )$$
and the estimated value of $\nu$ from the model is $\nu=0.48149$, as it should be, since we built the results assuming the value in~\eqref{nufromdata}. In rating chess players, one usually uses not $\bm{\sigma}$ itself, but ${\bf r}=(r_B,r_C,r_D)=400\log_{10}(\bm{\sigma})$, which we will again normalize by subtracting $r_C$ from ${\bf r}$ to make $r_C=0$ (this just corresponds to a different normalization than $\sum_i \sigma_i=1$). The resulting ratings are
$$ {\bf r}=(238.0135,0,140.3728)\,,$$
implying that de la Bourdonnais had a rating $238$ points above Cochrane's, and that Deschapelles, even playing with the handicap, had a rating $140$ points above Cochrane's, while de la Bourdonnais was $98$ points above the handicapped Deschapelles.

\section*{Acknowledgements}
This research did not receive any specific grant from funding agencies in the public, commercial, or not-for-profit sectors.


\end{document}